\documentclass[final,3p]{elsarticle}

\usepackage{amsmath,amssymb,amsfonts,mathtools,mathrsfs}
\usepackage{graphicx,float,subfigure}
\usepackage{xcolor}
\usepackage{algorithm,algorithmic}
\usepackage{pgfplots}
\usepackage{quiver}
\usepackage{lineno}
\usepackage{theorem}
\usepackage[colorlinks=true,breaklinks=true,linkcolor=black]{hyperref}

\usetikzlibrary{arrows}
\pgfplotsset{compat=1.18}
\modulolinenumbers[1]

\journal{CAGD}
\biboptions{authoryear}

\newtheorem{theorem}{Theorem}[section]
\newtheorem{proposition}{Proposition}[section]
\newtheorem{corollary}{Corollary}[section]
\newtheorem{lemma}{Lemma}[section]
\newtheorem{definition}{Definition}[section]
\newtheorem{remark}{Remark}[section]
\newtheorem{example}{Example}[section]
\newtheorem{conjecture}{Conjecture}[section]
\newenvironment{proof}{\par\noindent\textit{Proof.}\ }{\par}

\begin{document}

\begin{frontmatter}

\title{\Large\bf A Preliminary Study on the Dimensional Stability Classification of Polynomial Spline Spaces over T-meshes }

\author[1]{Bingru Huang}\corref{cor1}
\ead{hbr999@ustc.edu.cn}
\author[1]{Falai Chen}
\ead{chenfl@ustc.edu.cn}
\address[1]{School of Mathematical Sciences, University of Science and Technology of China, Hefei, 230026, People's Republic of China}
\cortext[cor1]{Corresponding author}

\begin{abstract}
This paper studies dimensional stability of polynomial spline spaces over T-meshes with highest-order smoothness. Dimensional stability is defined as the invariance of the dimension of the spline space over structurally isomorphic T-meshes, where \textcolor{black}{a structurally isomorphic map} preserves edge intersections and the mutual positions of all edges.
Absolute stability, a stronger notion, is introduced via structurally similar maps that depend solely on the topology of T-connected components.
It is proved that every T-connected component admits a unique decomposition into a diagonalizable part and a completely non-diagonalizable component (CNDC).
This reduces stability to the rank stability of the conformality matrix associated with the CNDC.
For diagonalizable T-meshes, the conformality vector space decomposes as a direct sum over independent T $l$-edges, supporting basis function construction per T $l$-edge.
These results provide a systematic framework for classifying dimensional stability and \textcolor{black}{a rigorous theoretical foundation for constructing} basis functions of spline spaces over T-meshes.
\end{abstract}


\end{frontmatter}

\section{Introduction}
Non-Uniform Rational B-Splines (NURBS) are the industry standard in Computer-Aided Design (CAD) for representing curves and surfaces. They offer high flexibility and precision for both analytic and free-form shapes using weighted control points and knot vectors. However, NURBS rely on tensor-product structures, requiring uniform global refinement for local changes. This increases computational cost and limits efficiency in applications needing adaptive local refinement. Spline formulations with local refinement capabilities \cite{ts1,HB1,thb,PHT,lr} overcome these limitations by allowing targeted refinement at specific cells. Such adaptive splines preserve spline continuity and smoothness while enhancing flexibility and efficiency through local refinement, reducing unnecessary computations.

Polynomial splines over T-meshes, a prominent class of adaptive splines, support local refinement at T-junctions, relaxing the rectangular grid requirement. A major challenge is determining the dimension of the spline space over arbitrary T-meshes, which is essential for basis construction and completeness in isogeometric analysis (IGA). For example, dimension formulas were used by~\cite{deng2008} to construct basis functions for hierarchical T-meshes (PHT-splines), demonstrating its critical role in developing practical adaptive spline bases. Several methodologies have been established for determining the dimension of spline spaces over T-meshes. These include the B-net method \cite{dim2006}, the Smoothing Cofactor Method \cite{dim20061,dim2016,zeng2015,zeng2016,zeng2018,bracco2019tchebycheffian,manni2018dimension,huang2023}, the Minimal Determining Set Method \cite{MDS,schumaker2012approximation}, the Homology Method \cite{dim2014,bracco2016dimension,bracco2016generalized,toshniwal2021polynomial,toshniwal2020dimension,toshniwal2021counting,toshniwal2023algebraic}, and the Space Embedding Method \cite{jin2013}, among others.

For spline spaces with low-order smoothness (smoothness order less than half the spline degree), explicit dimension formulas can be obtained by all the aforementioned methods. However, for high-order smoothness, no universal dimension formula exists, increasing the challenge. Recent work using the smoothing cofactor method introduces the concept of the completely non-diagonalizable component (CNDC), which reduces the dimension computation to determining the rank of the conformality matrix associated with the CNDC~\cite{huang2023}. The homology method offers a general dimension formula using homological techniques, but one component of the formula remains computationally intensive.

Research by~\cite{Ins2011} revealed dimensional instability in spline spaces over certain T-meshes. Dimensional instability means that the dimension depends not only on the topological structure of the T-mesh but also on its specific geometric configuration (e.g., vertex positions). This instability was further explored in \cite{Ins2012,li2019instability,guo2015problem}. These findings prompted investigations into T-meshes with stable dimensions, such as weighted T-meshes \cite{dim2014} and diagonalizable T-meshes \cite{dim2016}. Several adaptive spline frameworks, such as PHT-splines \cite{PHT}, HB-splines \cite{HB1,HB2}, THB-splines \cite{thb}, and LR-splines \cite{lr}, are dimensionally stable under their respective refinement rules.

This characteristic raises a central research question: Can T-meshes be systematically classified according to their topology to characterize those for which the spline space dimension is stable (depending solely on topology and independent of geometric configuration), thereby enabling the construction of a basis for the spline space over any such stable T-mesh?

{\color{black}Throughout this paper, an $l$-edge is a maximal connected union of collinear mesh edges. A T $l$-edge is an interior $l$-edge whose two endpoints are interior vertices. These notions are defined precisely in Section~\ref{sec.prelim}.}

This paper addresses the dimensional stability of polynomial spline spaces over T-meshes with highest-order smoothness through three main lines of work. We define dimensional stability as the invariance of the spline space dimension within the class of structurally isomorphic T-meshes, where a structural isomorphism preserves edge intersections and the mutual positions of parallel $l$-edges, and show that it depends on the rank stability of the conformality matrix. We then introduce absolute stability, a stricter variant based on structurally similar maps that preserve the \textcolor{black}{incidence structure and the number of vertices on each $l$-edge of a generalized T-connected component}, while allowing orientation swaps. Finally, we prove that every T-connected component admits a unique decomposition into diagonalizable and completely non-diagonalizable components, thereby reducing the stability problem to the rank stability of the conformality matrix associated with CNDC. For diagonalizable T-meshes, every T-connected component can be decomposed independently into individual T $l$-edges, thereby mapping the construction of spline basis functions directly onto each individual T $l$-edge. These developments establish a systematic framework for classifying dimensional stability and support the construction of stable bases for spline spaces over T-meshes.

The main results and advancements are:

\begin{itemize}
    \item We establish a formal mathematical framework for characterizing dimensional stability.
    \item We introduce a stronger notion termed absolute stability, which depends solely on the topological structure of a T-mesh.
    \item We prove the uniqueness of the complete partition proposed in \cite{huang2023}, which decomposes a T-connected component into a stable diagonalizable component and a completely non-diagonalizable component.
    \item We demonstrate that for diagonalizable T-meshes \cite{dim2016}, any T-connected component can be independently decomposed along each T $l$-edge, providing a rigorous theoretical foundation for the independent construction of spline basis functions.
\end{itemize}

The remainder of this paper is organized as follows. Section 2 reviews key concepts related to T-meshes, spline spaces over T-meshes, and related definitions. Section 3 introduces a structurally isomorphic map to define dimensional stability as an invariant within its class, depending on the rank stability of the associated conformality matrix, and proposes absolute stability using structurally similar maps to provide an initial classification of T-mesh dimensional stability. Section 4 examines the decomposition of T-connected components, introducing the $k$-partition and proving the existence and uniqueness of the complete partition. It also connects dimensional stability to the rank stability of the conformality matrix for CNDC and shows how diagonalizable T-meshes can be decomposed into independent T $l$-edges. Section 5 concludes the paper with a summary and directions for future research.
\section{Preliminaries}\label{sec.prelim}
\setcounter{equation}{0}
\setcounter{definition}{0}
In this section, we briefly introduce the concept of T-meshes and spline spaces over T-meshes, along with the smoothing cofactor method—a classic algebraic approach for determining the dimension of spline spaces.

\subsection{T-mesh, spline space over T-mesh and smoothing cofactor method}

\begin{definition}[T-mesh \cite{toshniwal2020,bracco2016dimension}]
    A T-mesh $\mathscr{T}=(\mathcal{T}_0,\mathcal{T}_1,\mathcal{T}_2)$ of $\mathbb{R}^2$ is defined as follows:
    \begin{itemize}
        \item[1.] $\mathcal{T}_2$ is a finite set of closed axis-aligned rectangles $f$ in $\mathbb{R}^2$, called cells;

        \item[2.] $\mathcal{T}_1 = \mathcal{T}_1^h \cup \mathcal{T}_1^v$ is a finite set of closed axis-aligned horizontal and vertical segments contained in $\bigcup_{f\in\mathcal{T}_2} \partial f$, called edges;

        \item[3.] $\mathcal{T}_0 \textcolor{black}{=} \bigcup_{l\in\mathcal{T}_1} \partial l$ is a finite set of points, called vertices.
    \end{itemize}

    Here, $\partial f$ denotes the boundary of the cell $f$, and $\partial l$ denotes the boundary of the edge $l$.

    These sets satisfy the following conditions:
    \begin{itemize}
        \item For each $f\in\mathcal{T}_2$, $\partial f$ is a finite union of elements of $\mathcal{T}_1$;

        \item For distinct $f,f'\in\mathcal{T}_2$, $f\cap f' = \partial f \cap \partial f'$ is a finite union of elements of $\mathcal{T}_0 \cup \mathcal{T}_1$;

        \item For distinct $l,l'\in\mathcal{T}_1$, $l\cap l' = \partial l \cap \partial l' \subseteq \mathcal{T}_0$;

        \item For each $v\in\mathcal{T}_0$, $v = l_1 \cap l_2$ where $l_1$ is a horizontal edge and $l_2$ is a vertical edge.
    \end{itemize}

    The domain of the T-mesh $\mathscr{T}$ is the region $\Omega \textcolor{black}{=} \bigcup_{f\in\mathcal{T}_2} f$.
\end{definition}

{\color{black}
Throughout this paper, the domain of $\mathscr{T}$ is assumed to be simply connected and its interior $\Omega^{0}$ is assumed to be connected. The vertices in $\Omega^{0}$ form the set $\mathcal{T}_0^{0}$ of interior vertices; the vertices in $\mathcal{T}_0\setminus\mathcal{T}_0^{0}$ are boundary vertices. Similarly, $\mathcal{T}_1^{0}$ denotes the set of interior edges, namely, the edges not contained in $\partial\Omega$.

A \emph{segment} of $\mathscr{T}$ is a connected union of collinear mesh edges. An $l$-edge (large edge) is a maximal such segment. An $l$-edge contained in the interior-edge set $\mathcal{T}_1^{0}$ is an \emph{interior $l$-edge}, whereas a maximal segment contained in $\partial\Omega$ is a \emph{boundary $l$-edge}. The set of boundary $l$-edges is denoted by $B(\mathscr{T})$. An interior $l$-edge is a \emph{cross-cut} if both endpoints are boundary vertices, a \emph{ray} if exactly one endpoint is a boundary vertex, and a \emph{T $l$-edge} if both endpoints are interior vertices. The sets of cross-cuts and rays are denoted by $C(\mathscr{T})$ and $R(\mathscr{T})$, respectively.

The edge--vertex configuration formed by all T $l$-edges and all vertices lying on them is denoted by $T(\mathscr{T})$. We assume that this configuration is connected and call it the \emph{T-connected component}; if it is disconnected, each connected component can be treated separately. A \emph{sub-component} of $T(\mathscr{T})$ is an edge--vertex substructure obtained by selecting T $l$-edges and specifying the vertices retained on those edges, with incidence inherited from $T(\mathscr{T})$. This convention permits a vertex shared by different selected edges to be assigned to an earlier component of a partition and omitted from a later one, as used in Section~\ref{sec: complete partition}.

A vertex lying on two T $l$-edges is a \emph{multi-vertex}; every other vertex on a T $l$-edge is a \emph{mono-vertex}. For an $l$-edge $l$, $n(l)$ denotes its number of retained vertices and $m(l)$ denotes its number of mono-vertices in the component under consideration.}

In Fig.~\ref{fig T-mesh}, the left subfigure illustrates a T-mesh denoted as $\mathscr{T}$. Here, among the vertices marked, $v_1, v_3$ and $v_4$ are boundary vertices, while the rest are interior vertices. Specifically, \textcolor{black}{$v_9$ and $v_{10}$ are} identified as multi-vertices, with all other interior vertices being mono-vertices. The line segment $v_1v_3$ forms a cross-cut, $v_4v_7$ forms a ray, whereas $v_2v_9$, $v_6v_{12}$, and $v_8v_{11}$ are T $l$-edges. The right subfigure shows the T-connected component of $\mathscr{T}$.

\begin{figure}
	\centering
	\begin{tikzpicture}[line cap=round,line join=round,>=triangle 45,x=1cm,y=1cm,scale=2]
		\draw [line width=1pt] (1,1)-- (3,1);
		\draw [line width=1pt] (3,1)-- (3,3.5);
		\draw [line width=1pt] (3,3.5)-- (1,3.5);
		\draw [line width=1pt] (1,3.5)-- (1,1);
		\draw [line width=1pt] (1,3)-- (3,3);
		\draw [line width=1pt] (1.5,3.5)-- (1.5,1);
		\draw [line width=1pt] (2.15,1.5)-- (2.15,2.5);
		\draw [line width=1pt] (1,2.5)-- (2.5,2.5);
		\draw [line width=1pt] (2.5,3.5)-- (2.5,1.5);
		\draw [line width=1pt] (1,1.5)-- (3,1.5);
		\draw [line width=1pt] (1.5,2)-- (2.5,2);
		\draw [line width=1pt] (1.86,3)-- (1.86,2);
		\draw [line width=1pt] (4,2)-- (5,2);
		\draw [line width=1pt] (4.36,3)-- (4.36,2);
        \draw [line width=1pt] (4.65,1.5)-- (4.65,2.5);
		\begin{scriptsize}
			\draw [fill=black] (1,3) circle (1pt);
			\draw[color=black] (1.1030303030302986,3.122626262626248) node {$v_1$};
			\draw [fill=black] (3,3) circle (1pt);
			\draw[color=black] (3.098989898989893,3.122626262626248) node {$v_3$};
			\draw [fill=black] (2.15,2.5) circle (1pt);
			\draw[color=black] (2.15,2.63) node {$v_6$};
            \draw [fill=black] (2.5,2.5) circle (1pt);
			\draw[color=black] (2.6,2.63) node {$v_7$};
            \draw [fill=black] (1,2.5) circle (1pt);
			\draw[color=black] (1.1030303030302986,2.63) node {$v_4$};
			\draw [fill=black] (1.5,2) circle (1pt);
			\draw[color=black] (1.6040404040403993,2.1286868686868577) node {$v_8$};
			\draw [fill=black] (2.5,2) circle (1pt);
			\draw[color=black] (2.6,2.1286868686868577) node {$v_{11}$};
            \draw [fill=black] (2.15,2) circle (1pt);
			\draw[color=black] (2.25,2.1286868686868577) node {$v_{10}$};
			\draw [fill=black] (1.86,3) circle (1pt);
			\draw[color=black] (1.86,3.122626262626248) node {$v_2$};
			\draw [fill=black] (1.86,2) circle (1pt);
			\draw[color=black] (1.98,2.1286868686868577) node {$v_9$};
			\draw [fill=black] (1.86,2.5) circle (1pt);
			\draw[color=black] (1.78,2.6296969696969566) node {$v_5$};
            \draw [fill=black] (2.15,1.5) circle (1pt);
			\draw[color=black] (2.28,1.6296969696969566) node {$v_{12}$};
			\draw [fill=black] (4.36,3) circle (1pt);
			\draw[color=black] (4.36,3.122626262626248) node {$v_2$};
			\draw [fill=black] (4.36,2.5) circle (1pt);
			\draw[color=black] (4.26,2.6296969696969566) node {$v_5$};
            \draw [fill=black] (4.36,2) circle (1pt);
			\draw[color=black] (4.48,2.1286868686868577) node {$v_9$};
            \draw [fill=black] (4,2) circle (1pt);
			\draw[color=black] (4.1040404040403993,2.1286868686868577) node {$v_8$};
            \draw [fill=black] (4.65,2.5) circle (1pt);
			\draw[color=black] (4.65,2.63) node {$v_{6}$};
            \draw [fill=black] (4.65,1.5) circle (1pt);
			\draw[color=black] (4.78,1.6296969696969566) node {$v_{12}$};
            \draw [fill=black] (5,2) circle (1pt);
			\draw[color=black] (5.1,2.1286868686868577) node {$v_{11}$};
            \draw [fill=black] (4.65,2) circle (1pt);
			\draw[color=black] (4.75,2.1286868686868577) node {$v_{10}$};
		\end{scriptsize}
	\end{tikzpicture}
	\caption{\label{fig T-mesh} T-mesh and T-connected component}
\end{figure}

\begin{definition}[Spline Space over T-meshes]
    Given a T-mesh $\mathscr{T}$, let $\Omega$ be the region occupied by the cells in $\mathscr{T}$, $\mathbf{d}=(d_1,d_2)$, $\mathbf{r}=(r_1,r_2)$. Then the spline space over $\mathscr{T}$ is defined as follows.
	$$\textcolor{black}{\mathbf{S}}_{\mathbf{d}}^{\mathbf{r}}(\mathscr{T})=\{s(x,y)\in C^{\mathbf{r}}(\Omega) \big|\quad s(x,y)|_{\phi}\in \mathbb{P}_{\mathbf{d}}, \forall\phi\in\mathcal{T}_2\}$$
	where $\mathbb{P}_{\mathbf{d}}$ is the function space of all polynomials with bi-degree $(d_1,d_2)$, $C^{\mathbf{r}}(\Omega)$ is the space consisting of all the bivariate functions continuous in $\Omega$ with order $r_1$ along the $x$-direction and with order $r_2$ along the $y$-direction. If $(r_1,r_2)=(d_1-1,d_2-1)$, then the spline space is called \textbf{highest-order smoothness spline space.} Assuming the same degree $d$ in both directions (i.e., $d_1=d_2=d$), this spline space is denoted by the simplified symbol $\mathbf{S}_{d}(\mathscr{T})$, \textcolor{black}{that is, $\mathbf{S}_{d}(\mathscr{T})=\mathbf{S}_{(d,d)}^{(d-1,d-1)}(\mathscr{T})$}.
\end{definition}

    In this paper, we mainly discuss the \textbf{\textcolor{black}{highest-order}} smoothness spline space $\mathbf{S}_{d}(\mathscr{T})$. In the following, the concept of diagonalizable T-meshes is introduced. First, notation is defined for an ordered sequence of all T $l$-edges, denoted $l_1 \succ l_2 \succ \dots \succ l_t$. The term $\overline{l}_i$ denotes the edge $l_i$ with intersection vertices shared with prior edges $l_j$ ($j = 1, 2, \dots, i-1$) removed, for $2 \le i \le t$, where $\overline{l}_1 = l_1$. This defines the set $\{\overline{l}_1, \overline{l}_2, \dots, \overline{l}_t\}$, termed the $\textbf{$t$-partition}$ of $T(\mathscr{T})$.

\begin{definition}[\cite{dim2016}]\label{diagonalizable}
For a given T-mesh $\mathscr{T}$, consider the spline space $\mathbf{S}_{d}(\mathscr{T})$. If there exists a $t$-partition of $T(\mathscr{T})$, denoted $\{\overline{l}_1, \overline{l}_2, \dots, \overline{l}_t\}$, such that $n(\overline{l}_i)\geq d+1$ for each $i=1,2,\dots,t$ (where $n(\overline{l}_i)$ denotes the number of vertices lying on the line segment $\overline{l}_i$), then $\mathscr{T}$ is called a $\textbf{diagonalizable T-mesh}$.
\end{definition}

Concerning the dimension of spline spaces over diagonalizable T-meshes, we present the following result from~\cite{dim2016}.

\begin{proposition}[\cite{dim2016}]
\label{dimension diagonalizable}
If a T-mesh $\mathscr{T}$ is diagonalizable, then the dimension of the spline space $\mathbf{S}_{d}(\mathscr{T})$ is
\begin{equation}\label{eq1}
	\dim \mathbf{S}_{d}(\mathscr{T})=(d+1)^2+(c-t)(d+1)+n_v,
\end{equation}
where $c$ is the number of cross-cuts of $\mathscr{T}$, $t$ is the number of T $l$-edges of $T(\mathscr{T})$ and $n_v$ is the number of interior vertices.
\end{proposition}

Proposition~\ref{dimension diagonalizable} demonstrates that the dimension of the spline space over a diagonalizable T-mesh is stable, indicating that the dimension formula~\eqref{eq1} depends solely on the topological characteristics of the T-mesh (like $c$, $t$, $n_v$) and not on its geometric position. We refer to such T-meshes, where the dimension of the spline space is stable, as \textbf{dimensionally stable T-meshes}.

The smoothing cofactor method is a classic algebraic approach to determine the dimension of spline spaces. It was first introduced in \cite{sm1,CHUI1983197,ChuiWang1983} and has been applied to T-meshes in \cite{dim20061,dim2016,wu2013,zeng2015,zeng2016,zeng2018,huang2023}. We briefly review the method below, with emphasis on the conformality conditions and rank analysis relevant to dimensional stability.

Let $s(x,y)\in \mathbf{S}_{d}(\mathscr{T})$ be a spline function. Referring to the left subfigure of Fig.~\ref{fig local conformaity condition}, let $f_i$ and $f_j$ be two adjacent cells in $\mathscr{T}$ with $i\neq j$. Suppose $s(x,y)|_{f_i}=s_i(x,y)$, $s(x,y)|_{f_j}=s_j(x,y)$, and their common edge $l_{ij}$ lies on the line $x=x_{i,j}$. By B\'ezout's theorem, there exists \textcolor{black}{a univariate polynomial $\gamma_{ij}(y)$ of degree at most $d$} such that
\[
s_i(x,y)-s_j(x,y)=\gamma_{ij}(y)(x-x_{i,j})^{d}.
\]
This shows that if one of $s_i(x,y)$ or $s_j(x,y)$ is known, the other is determined by $\gamma_{ij}(y)$, which has fewer free variables. We call $\gamma_{ij}(y)$ an \textbf{edge cofactor} of $s(x,y)$ along $l_{ij}$.

\begin{figure}
\centering
\definecolor{uuuuuu}{rgb}{0.26666666666666666,0.26666666666666666,0.26666666666666666}
\begin{tikzpicture}[line cap=round,line join=round,>=triangle 45,x=1cm,y=1cm,scale=0.8]
\clip(-2.73939393939394,-1.4036363636363638) rectangle (9.616161616161603,6.7094949494949265);
\draw [line width=1pt] (0,4.5)-- (0,0.5);
\draw (-1.6404040404040423,3.0408080808080684) node[anchor=north west] {$f_i$};
\draw (1.041414141414137,3.0408080808080684) node[anchor=north west] {$f_j$};
\draw (-2.0505050505050522,2.0064646464646367) node[anchor=north west] {$s_i(x,y)$};
\draw (0.641414141414137,2.0064646464646367) node[anchor=north west] {$s_j(x,y)$};
\draw (-0.64707070707071005,0.4468686868686816) node[anchor=north west] {$x=x_{i,j}$};
\draw [line width=1pt] (4,2.5)-- (8,2.5);
\draw [line width=1pt] (6,4.5)-- (6,0.5);
\draw (4.824242424242416,3.7761616161616014) node[anchor=north west] {$f_1$};
\draw (4.824242424242416,1.764040404040395) node[anchor=north west] {$f_2$};
\draw (6.503030303030293,1.7397979797979708) node[anchor=north west] {$f_3$};
\draw (6.503030303030293,3.7761616161616014) node[anchor=north west] {$f_4$};
\draw (8.110101010100999,2.7822222222222104) node[anchor=north west] {$y=y_0$};
\draw (5.4,4.980202020202002) node[anchor=north west] {$x=x_0$};
\draw (6.004040404040396,2.4610101010100897) node[anchor=north west] {$v$};
\begin{scriptsize}
\draw [fill=uuuuuu] (6,2.5) circle (2pt);
\end{scriptsize}
\end{tikzpicture}
\caption{\label{fig local conformaity condition} Local conformality condition}
\end{figure}

From a local perspective, an interior vertex is the intersection of two $l$-edges. Referring to the right subfigure of Fig.~\ref{fig local conformaity condition}, let an interior vertex $v$ be intersected by an $x$-edge $l_1$ ($y=y_0$) and a $y$-edge $l_2$ ($x=x_0$), with $s(x,y)|_{f_i}=s_i(x,y)$ for $i=1,2,3,4$. Then
\begin{gather}\label{eq2}
\gamma_{43}(x)-\gamma_{12}(x)=\delta(x-x_0)^{d},\\ \label{eq3}
\gamma_{41}(y)-\gamma_{32}(y)=\delta(y-y_0)^{d}.
\end{gather}
If one of $\gamma_{12}(x)$ or $\gamma_{43}(x)$ is given, the other is determined by $\delta$, which has fewer free variables (equation (\ref{eq3}) is analogous). We call $\delta$ a \textbf{vertex cofactor} of $s(x,y)$ corresponding to $v$.

From the global perspective, interior vertex cofactors are not mutually linearly independent. For an $l$-edge of $\mathscr{T}$, its edge cofactors form a $C^{d-1}$ continuous univariate spline space \cite{wu2013}. For a T $l$-edge, vertex constraints apply. For example, on a horizontal T $l$-edge with $r$ vertices at $x_1,\dots,x_r$ and vertex cofactors $\delta_1,\dots,\delta_r$ (Fig.~\ref{fig gcc}), let the edge cofactors from left to right be $\gamma_1, \dots, \gamma_{r-1}$. According to the definition of vertex cofactors as the relations between adjacent edge cofactors, we can establish the following system of equations:
\begin{equation*}
\left\{
\begin{aligned}
\gamma_1 - 0 &= \delta_1(x-x_1)^d, \\
\gamma_2 - \gamma_1 &= \delta_2(x-x_2)^d, \\
&\ \ \vdots \\
\gamma_{r-1} - \gamma_{r-2} &= \delta_{r-1}(x-x_{r-1})^d, \\
0 - \gamma_{r-1} &= \delta_r(x-x_r)^d.
\end{aligned}
\right.
\end{equation*}
Summing all these equations together naturally cancels out all the intermediate edge cofactors $\gamma_i$, which directly yields
\begin{equation}\label{eq4}
\sum_{i=1}^{r} \delta_i (x-x_i)^d = 0.
\end{equation}

\begin{figure}
    \centering
\begin{tikzpicture}[line cap=round,line join=round,>=triangle 45,x=1cm,y=1cm,scale=2]
\draw [line width=1pt] (1,0.5)-- (1,2.5);
\draw [line width=1pt] (1,1.5)-- (5,1.5);
\draw [line width=1pt] (5,2.5)-- (5,0.5);

\draw [line width=1pt] (1.5,1.5)-- (1.5,2);
\draw [line width=1pt] (2,1.5)-- (2,1);
\draw [line width=1pt] (2.5,2)-- (2.5,1);
\draw [line width=1pt] (3,1.5)-- (3,1);
\draw [line width=1pt] (4.5,1.5)-- (4.5,1);
\draw [line width=1pt] (4,1.5)-- (4,2);

\draw (1.25,1.5) node[anchor=south] {$\gamma_1$};
\draw (1.75,1.5) node[anchor=south] {$\gamma_2$};
\draw (2.25,1.5) node[anchor=south] {$\gamma_3$};
\draw (2.75,1.5) node[anchor=south] {$\gamma_4$};
\draw (3.5,1.5)  node[anchor=south] {$\ldots$}; 
\draw (4.25,1.5) node[anchor=south] {$\gamma_{r-2}$};
\draw (4.75,1.5) node[anchor=south] {$\gamma_{r-1}$};

\draw (0.94,1.5) node[anchor=north west] {$\delta_1$};
\draw (1.35,1.5) node[anchor=north west] {$\delta_2$};
\draw (1.95,1.5) node[anchor=north west] {$\delta_3$};
\draw (2.48,1.5) node[anchor=north west] {$\delta_4$};
\draw (2.94,1.5) node[anchor=north west] {$\delta_5$};
\draw (3.82,1.5) node[anchor=north west] {$\delta_{r-2}$};
\draw (4.47,1.5) node[anchor=north west] {$\delta_{r-1}$};
\draw (4.95,1.5) node[anchor=north west] {$\delta_r$};
\draw (3.32,1.5) node[anchor=north west] {$\ldots$};

\begin{scriptsize}
\draw [fill=black] (1,1.5) circle (1pt);
\draw [fill=black] (1.5,1.5) circle (1pt);
\draw [fill=black] (2,1.5) circle (1pt);
\draw [fill=black] (2.5,1.5) circle (1pt);
\draw [fill=black] (3,1.5) circle (1pt);
\draw [fill=black] (4,1.5) circle (1pt);
\draw [fill=black] (4.5,1.5) circle (1pt);
\draw [fill=black] (5,1.5) circle (1pt);
\end{scriptsize}
\end{tikzpicture}
    \caption{\label{fig gcc}Vertex cofactors and edge cofactors along a horizontal T $l$-edge}
\end{figure}

This equation is equivalent to the linear system (denoted $\mathscr{P}=0$):
\begin{equation}
\begin{pmatrix}
1 & 1 & \cdots & \cdots & 1 \\
x_1 & x_2 & \cdots & \cdots & x_r\\
x_1^2 & x_2^2 & \cdots & \cdots & x_r^2\\
\vdots & \vdots & \ddots & \ddots & \vdots\\
x_1^{d-1} & x_2^{d-1} & \cdots & \cdots & x_r^{d-1}\\
x_1^d & x_2^d & \cdots & \cdots & x_r^d
\end{pmatrix}
\begin{pmatrix}
\delta_1 \\
\delta_2 \\
\vdots \\
\delta_{r-1} \\
\delta_{r}
\end{pmatrix}
=\begin{pmatrix}
0 \\
\vdots \\
0
\end{pmatrix}.
\label{gcc}
\end{equation}

The linear systems for all T $l$-edges form the \textbf{global conformality condition} of $\mathbf{S}_{d}(\mathscr{T})$. When $r \le d+1$, the null space of \eqref{gcc} is trivial; such a T $l$-edge contributes nothing to the dimension and is called a \textbf{vanishable edge}. We assume no vanishable T $l$-edges exist in the T-mesh.

The conformality vector space (CVS) of a T-connected component is defined as follows \cite{zeng2015}.
\begin{definition}\label{def2.4}
\textcolor{black}{Let $T(\mathscr{T})$ be the T-connected component of $\mathscr{T}$. Denote its T $l$-edges by $l_1,\ldots,l_t$ and its vertex cofactors by $\delta_1,\ldots,\delta_v$. The conformality vector space (\textsf{CVS}) of $T(\mathscr{T})$ is}
\[
\textsf{CVS}[T(\mathscr{T})]\textcolor{black}{=}
\left\{
\boldsymbol{\delta}=(\delta_1,\ldots,\delta_v)
\,\middle|\,
\mathscr{P}_{l_i}=0\ \text{for }1\le i\le t
\right\}.
\]
\textcolor{black}{Here $\mathscr{P}_{l_i}=0$ denotes the conformality condition across the T $l$-edge $l_i$, analogous to equation~\eqref{gcc}. The coefficient matrix of these conditions is called the \textbf{conformality matrix} and is denoted by $M(T(\mathscr{T}))$.}
\end{definition}

The dimension formula follows \cite{zeng2015}.
\begin{theorem}\label{thm2.1}
Given a T-mesh $\mathscr{T}$, with $M(T(\mathscr{T}))$ the conformality matrix of $T(\mathscr{T})$, we have
\begin{equation}\label{eqdim}
    \dim \mathbf{S}_{d}(\mathscr{T})=(d+1)^2+c(d+1)+n_v-\text{rank}(M(T(\mathscr{T}))),
\end{equation}
where $c$ is the number of cross-cuts and $n_v$ the number of interior vertices.
\end{theorem}

Thus, studying $\dim \mathbf{S}_{d}(\mathscr{T})$ reduces to analyzing $\text{rank}(M(T(\mathscr{T})))$. The concept of complete partition separates the diagonalizable and non-diagonalizable parts of a T-connected component; an algorithm computes the complete partition \cite{huang2023}. These are central to understanding dimensional instability here.

\begin{definition}[\cite{huang2023}]\label{def bipar}
    Let $T(\mathscr{T})$ be the T-connected component of a T-mesh $\mathscr{T}$ and let \textcolor{black}{$S\subset T(\mathscr{T})$} be a sub-component. Then we call \textcolor{black}{$\{S,T(\mathscr{T})\setminus S\}$} a \textcolor{black}{bipartition} of $T(\mathscr{T})$.
\end{definition}

\begin{definition}[\cite{huang2023}]
Let $T(\mathscr{T})$ be the T-connected component of a T-mesh $\mathscr{T}$, and $\{\overline{T_1},\overline{T_2}\}$ a bipartition with $|\overline{T_1}|=s$. If $\overline{T_1}$ has no $s$-partition with $n(\overline{l})>d+1$ for all $\overline{l}\in\overline{T_1}$ while $\overline{T_2}$ does, then $\overline{T_1}$ is a \textbf{non-diagonalizable component}, $\overline{T_2}$ a \textbf{diagonalizable component}, and $\{\overline{T_1},\overline{T_2}\}$ a \textbf{regular partition} of $T(\mathscr{T})$.
\end{definition}

\begin{definition}[\cite{huang2023}]
Let $\mathscr{T}$ be a T-mesh with T-connected component $T(\mathscr{T})$. If $\{\overline{T_1},\overline{T_2}\}$ is a regular partition and no other regular partition $\{\overline{T'_1},\overline{T'_2}\}$ satisfies $\overline{T'_1}\subseteq \overline{T_1}$, then $\{\overline{T_1},\overline{T_2}\}$ is a \textbf{complete partition} of $T(\mathscr{T})$, and $\overline{T_1}$ its \textbf{completely non-diagonalizable component} (CNDC).
\end{definition}

The following algorithm computes a complete partition \cite{huang2023}.

\begin{algorithm}[htb]
\caption{Decompose the T-connected component of a T-mesh into a complete partition.}
\label{alg}
\begin{algorithmic}[1]
\REQUIRE A T-connected component $T(\mathscr{T})$ of a given T-mesh $\mathscr{T}$;
\ENSURE The complete partition of $T(\mathscr{T})$;
\STATE $T_1=T(\mathscr{T})$, $T_2=\varnothing$;
\REPEAT
\FOR{$l\in T_1$}
\IF{ $m(l)\ge d+1$}
\STATE Set $T_1=T_1\backslash\{\tilde{l}\}$, where $\tilde{l}$ is $l$ with multi-vertices removed;
\STATE Set $T_2=T_2\cup\{\tilde{l}\}$;
\ENDIF
\ENDFOR
\UNTIL{$T_1=\varnothing$ or $m(l)<d+1$ for any $l\in T_1$}
\STATE Output $\{T_1,T_2\}$.
\end{algorithmic}
\end{algorithm}

Dimensional instability of spline spaces over T-meshes reduces to the rank of the conformality matrix associated with the CNDC \cite{huang2023}.

\subsection{Summary of Section}~\ref{sec.prelim}
In this section, we have reviewed the definition of T-meshes, polynomial spline spaces over T-meshes, and the smoothing cofactor method as a key algebraic tool for computing the dimension of such spaces. Building on earlier developments of this method, we have also presented the concept of diagonalizable T-meshes, which ensure dimension stability under highest-order smoothness, as well as the conformality vector space associated with T-connected components. Furthermore, recent work by Huang and Chen (2024) provides a refined framework by decomposing T-connected components into diagonalizable and CNDC through a complete partition process. This decomposition leads to a more precise dimension formula, ultimately reducing the computation of the spline space dimension to determining the rank—or equivalently, the dimension—of the conformality vector space over the CNDC.
Since a substantial number of new notations and concepts are introduced here, we summarize the main symbols used in this section in Table~\ref{tab:symbols}. Corresponding symbol tables will be updated or supplemented in subsequent sections as needed to maintain clarity and consistency throughout the paper.
\begin{table}[H]
\centering
\caption{List of Symbols in Section~\ref{sec.prelim}}
\begin{tabular}{|c|l|}
\hline
\textbf{Symbol} & \textbf{Explanation} \\
\hline
$\mathscr{T}$ & T-mesh \\
\hline
$f\in\mathcal{T}_2$ & A cell $f$ of $\mathscr{T}$ \\
\hline
$\textcolor{black}{S\subset T(\mathscr{T})}$ & $S$ is a sub-component of $T(\mathscr{T})$ \\
\hline
$l\in\mathcal{T}_1$ & An edge $l$ of $\mathscr{T}$ \\
\hline
$v\in\mathcal{T}_0$ & A vertex $v$ of $\mathscr{T}$\\
\hline
$\Omega$ & The domain occupied by $\mathscr{T}$\\
\hline
$B(\mathscr{T})$ & The set of all boundary $l$-edges of $\mathscr{T}$\\
\hline
$C(\mathscr{T})$ & The set of all cross-cuts of $\mathscr{T}$\\
\hline
$R(\mathscr{T})$ & The set of all rays of $\mathscr{T}$\\
\hline
$T(\mathscr{T})$ & \textcolor{black}{The T-connected component of $\mathscr{T}$}\\
\hline
$\textcolor{black}{\mathbf{S}}_{\mathbf{d}}^{\mathbf{r}}(\mathscr{T})$ & Spline space with bi-degree $\mathbf{d}$ and smoothness order $\mathbf{r}$ over $\mathscr{T}$\\
\hline
$\textcolor{black}{\mathbf{S}}_{d}(\mathscr{T})$ & Spline space with bi-degree \textcolor{black}{$(d,d)$} and highest-order smoothness\\
\hline
$l_i$ & The $i$-th $l$-edge \\
\hline
$\overline{l}_i$ & The $l$-edge $l_i$ after removing intersection points with previous $l$-edges \\
\hline
$n(l)$ & Number of vertices on $l$-edge $l$\\
\hline
$m(l)$ & The number of mono-vertices on $l$ \\
\hline
$c$ & Number of cross-cuts in the T-mesh \\
\hline
$t$ & Number of T $l$-edges in the T-mesh \\
\hline
$n_v$ & Number of interior vertices in the T-mesh \\
\hline
$\gamma$ & An edge cofactor\\
\hline
$\delta$ & A vertex cofactor\\
\hline
$\textsf{CVS}[T(\mathscr{T})]$ & The conformality vector space of $T(\mathscr{T})$\\
\hline
$M(T(\mathscr{T}))$ & The conformality matrix of $T(\mathscr{T})$\\
\hline
\end{tabular}
\label{tab:symbols}
\end{table}


\section{Structurally Similar Map and Absolute Stability of T-meshes}\label{sec:Structurally Similar Map}
This section defines dimensional stability for T-meshes.
In Subsection~\ref{subsection definition stab}, we define dimensional stability via the structurally isomorphic map, an equivalence relation on the set $U$ of all T-meshes.
In Subsection~\ref{subsection: defintion ab}, we introduce absolute stability, a stricter notion.

\subsection{The definition of stability of T-meshes}\label{subsection definition stab}
In mathematics, stability is usually described as the invariance of a quantity during a change. In this section, \textcolor{black}{we describe} the dimensional stability of a T-mesh as an invariant on a class of T-meshes. To this end, we first give a \textcolor{black}{structurally isomorphic map} between two T-meshes as follows.

\begin{definition}[Structurally Isomorphic Map]\label{defsim}
\color{black}
Let $\mathscr{T}_1$ and $\mathscr{T}_2$ be two T-meshes. Denote by $L(\mathscr{T}_i)$ the set of all $l$-edges of $\mathscr{T}_i$, including its boundary $l$-edges, and by $L_1(\mathscr{T}_i)$ and $L_2(\mathscr{T}_i)$ its horizontal and vertical $l$-edges, respectively. For an axis-aligned $l$-edge $l$, let $\operatorname{cor}(l)$ be its constant $y$-coordinate when $l$ is horizontal and its constant $x$-coordinate when $l$ is vertical.

The T-meshes $\mathscr{T}_1$ and $\mathscr{T}_2$ are \textbf{structurally isomorphic} if there exists a bijection
\[
\alpha:L(\mathscr{T}_1)\longrightarrow L(\mathscr{T}_2)
\]
satisfying the following conditions:
\begin{itemize}
    \item $\alpha$ preserves the types of special $l$-edges:
    \[
    \alpha(B(\mathscr{T}_1))=B(\mathscr{T}_2),\quad
    \alpha(C(\mathscr{T}_1))=C(\mathscr{T}_2),\quad
    \alpha(R(\mathscr{T}_1))=R(\mathscr{T}_2),\quad
    \alpha\bigl(L(T(\mathscr{T}_1))\bigr)=L(T(\mathscr{T}_2)).
    \]
    \item $\alpha$ preserves incidence: for all $l_i,l_j\in L(\mathscr{T}_1)$,
    \[
    l_i\cap l_j\ne\varnothing
    \quad\Longleftrightarrow\quad
    \alpha(l_i)\cap\alpha(l_j)\ne\varnothing.
    \]
    \item $\alpha$ preserves the mutual positions of parallel $l$-edges: for $k=1,2$ and $l_i,l_j\in L_k(\mathscr{T}_1)$,
    \[
    \operatorname{cor}(l_i)<\operatorname{cor}(l_j)
    \quad\Longrightarrow\quad
    \operatorname{cor}(\alpha(l_i))<\operatorname{cor}(\alpha(l_j)).
    \]
\end{itemize}
Such a bijection is called a \textbf{structurally isomorphic map} from $\mathscr{T}_1$ to $\mathscr{T}_2$.
\end{definition}

\textcolor{black}{By Definition~\ref{defsim}, we obtain the following lemma.}

\begin{lemma}~\label{lemma h to v}
Let $\alpha: L(\mathscr{T}_1) \to L(\mathscr{T}_2)$ be a structurally isomorphic map.
If there exists an $l$-edge $l_0$ of $\mathscr{T}_1$ that is horizontal (resp. vertical) such that $\alpha(l_0)$ is vertical (resp. horizontal) in $\mathscr{T}_2$,
then $\alpha$ maps every horizontal $l$-edge of $\mathscr{T}_1$ to a vertical $l$-edge of $\mathscr{T}_2$, and every vertical $l$-edge of $\mathscr{T}_1$ to a horizontal $l$-edge of $\mathscr{T}_2$. In other words, $\alpha$ either preserves the horizontal/vertical orientation of all $l$-edges or reverses it for all $l$-edges.
\end{lemma}

\begin{proof}
\color{black}
Assume, without loss of generality, that $l_0$ is horizontal and $\alpha(l_0)$ is vertical. Consider the incidence graph whose nodes are the $l$-edges of $\mathscr{T}_1$ and whose adjacent nodes correspond to intersecting $l$-edges. This graph is connected because the T-mesh is connected. Moreover, adjacent nodes always represent perpendicular $l$-edges.

Let $l$ be any $l$-edge of $\mathscr{T}_1$. Choose an incidence path
\[
l_0=e_0,e_1,\ldots,e_q=l.
\]
By incidence preservation, $\alpha(e_0),\alpha(e_1),\ldots,\alpha(e_q)$ is an incidence path in $\mathscr{T}_2$. Orientations alternate along both paths. Since the orientation of $\alpha(e_0)$ is opposite to that of $e_0$, induction along the path shows that the orientation of $\alpha(e_i)$ is opposite to that of $e_i$ for every $i$. In particular, every horizontal $l$-edge is mapped to a vertical one and every vertical $l$-edge is mapped to a horizontal one. If no $l$-edge changes orientation, the map preserves both orientation classes. Thus the map either preserves all orientations or reverses all of them.

$\Box$
\end{proof}

By Lemma~\ref{lemma h to v}, T-meshes structurally isomorphic to a given T-mesh $  \mathscr{T}  $ fall into exactly two categories. In the first category, the structurally isomorphic map preserves the horizontal/vertical orientation of all $  l  $-edges, so horizontal $  l  $-edges of $  \mathscr{T}  $ correspond bijectively to horizontal $  l  $-edges of the isomorphic mesh, and vertical to vertical. In the second category, the map reverses the orientation, so horizontal $  l  $-edges of $  \mathscr{T}  $ correspond bijectively to vertical $  l  $-edges of the isomorphic mesh, and vertical $  l  $-edges to horizontal.
The three conditions in Definition~\ref{defsim} ensure structural equivalence: the first preserves the types of special $  l  $-edges (boundary, cross-cuts, rays, and T $l$-edges); the second preserves intersection relations between $  l  $-edges; the third preserves the \textcolor{black}{mutual positions} within each direction (horizontal or vertical).

Together, these conditions guarantee that any T-mesh structurally isomorphic to $\mathscr{T}$ has the same $l$-edge connectivity and \textcolor{black}{mutual positions of parallel $l$-edges} as $\mathscr{T}$, differing only in the geometric positions of the corresponding edges. The two categories differ precisely by a $\pi/2$ rotation in the orientation of the $l$-edge connectivity.

\begin{proposition}\label{proof of equivalence of sim}
    \textcolor{black}{Let $U$ be the set of all T-meshes and define a relation $\sim_{\mathrm I}$ on $U$ by $\mathscr{T}_1\sim_{\mathrm I}\mathscr{T}_2$ if and only if $\mathscr{T}_1$ is structurally isomorphic to $\mathscr{T}_2$. Then $\sim_{\mathrm I}$ is an equivalence relation on $U$.}
\end{proposition}
\begin{proof}
\color{black}
Reflexivity follows from the identity map on $L(\mathscr{T})$. If $\alpha:L(\mathscr{T}_1)\to L(\mathscr{T}_2)$ is a structurally isomorphic map, then its inverse $\alpha^{-1}$ preserves edge types, incidence, and mutual positions; hence symmetry holds.

For transitivity, suppose that $\alpha_1:L(\mathscr{T}_1)\to L(\mathscr{T}_2)$ and $\alpha_2:L(\mathscr{T}_2)\to L(\mathscr{T}_3)$ are structurally isomorphic maps. Their composition $\alpha=\alpha_2\circ\alpha_1$ is a bijection. Each special edge type is preserved successively by $\alpha_1$ and $\alpha_2$. For any $l_i,l_j\in L(\mathscr{T}_1)$,
\[
l_i\cap l_j\ne\varnothing
\Longleftrightarrow
\alpha_1(l_i)\cap\alpha_1(l_j)\ne\varnothing
\Longleftrightarrow
\alpha(l_i)\cap\alpha(l_j)\ne\varnothing.
\]
Likewise, applying the mutual-position condition first to $\alpha_1$ and then to $\alpha_2$ shows that $\alpha$ preserves the mutual positions of parallel $l$-edges. Thus $\alpha$ is structurally isomorphic, proving transitivity. Therefore $\sim_{\mathrm I}$ is an equivalence relation.

$\Box$
\end{proof}

\begin{definition}
    Let $U$ be the set of all T-meshes. Then, any equivalence class induced by the structurally isomorphic map is called a \textbf{structurally isomorphic class}. The set of all structurally isomorphic classes is denoted by $\mathrm{SIC}(U)$ ($\mathrm{SIC}$ in short).
\end{definition}

\begin{example}
    \textcolor{black}{Consider the T-mesh $\mathscr{T}$ shown in Fig.~\ref{fig T-mesh}. Let $[\mathscr{T}]_{\mathrm{I}}\in \mathrm{SIC}$ be the class of all T-meshes structurally isomorphic to $\mathscr{T}$. Figure~\ref{fig structurally isomorphic class} shows two T-meshes $\mathscr{T}_1,\mathscr{T}_2\in[\mathscr{T}]_{\mathrm{I}}$.}
\end{example}
\begin{figure}[H]
    \centering
    \begin{tikzpicture}[line cap=round,line join=round,>=triangle 45,x=1.0cm,y=1.0cm,scale=0.7]
\draw [line width=1.pt] (1.,1.)-- (1.,7.);
\draw [line width=1.pt] (1.,7.)-- (12.,7.);
\draw [line width=1.pt] (12.,7.)-- (12.,1.);
\draw [line width=1.pt] (12.,1.)-- (1.,1.);
\draw [line width=1.pt] (1.,6.)-- (12.,6.);
\draw [line width=1.pt] (1.,2.)-- (12.,2.);
\draw [line width=1.pt] (3.,7.)-- (3.,1.);
\draw [line width=1.pt] (10.,7.)-- (10.,2.);
\draw [line width=1.pt] (1.,4.)-- (10.,4.);
\draw [line width=1.pt] (3.,3.)-- (10.,3.);
\draw [line width=1.pt] (6.,6.)-- (6.,3.);
\draw [line width=1.pt] (9.,2.)-- (9.,4.);
\draw [line width=1.pt] (16.,-1.)-- (16.,9.);
\draw [line width=1.pt] (16.,9.)-- (24.,9.);
\draw [line width=1.pt] (24.,9.)-- (24.,-1.);
\draw [line width=1.pt] (24.,-1.)-- (16.,-1.);
\draw [line width=1.pt] (18.,9.)-- (18.,-1.);
\draw [line width=1.pt] (23.,9.)-- (23.,-1.);
\draw [line width=1.pt] (16.,6.)-- (23.,6.);
\draw [line width=1.pt] (16.,1.)-- (24.,1.);
\draw [line width=1.pt] (19.,6.)-- (19.,-1.);
\draw [line width=1.pt] (21.,6.)-- (21.,1.);
\draw [line width=1.pt] (19.,4.)-- (23.,4.);
\draw [line width=1.pt] (18.,2.)-- (21.,2.);
\begin{scriptsize}
\draw [fill=black] (1.,6.) circle (3.0pt);
\draw[color=black] (1.25,6.3664) node {$v'_1$};
\draw [fill=black] (12.,6.) circle (3.0pt);
\draw[color=black] (12.25,6.3664) node {$v'_3$};
\draw [fill=black] (1.,4.) circle (3.0pt);
\draw[color=black] (1.25,4.3644) node {$v'_4$};
\draw [fill=black] (10.,4.) circle (3.0pt);
\draw[color=black] (10.25,4.3644) node {$v'_7$};
\draw [fill=black] (3.,3.) circle (3.0pt);
\draw[color=black] (3.25,3.3524) node {$v'_8$};
\draw [fill=black] (10.,3.) circle (3.0pt);
\draw[color=black] (10.3,3.3524) node {$v'_{11}$};
\draw [fill=black] (6.,6.) circle (3.0pt);
\draw[color=black] (6.25,6.3664) node {$v'_2$};
\draw [fill=black] (6.,3.) circle (3.0pt);
\draw[color=black] (6.25,3.3524) node {$v'_9$};
\draw [fill=black] (9.,2.) circle (3.0pt);
\draw[color=black] (9.3,2.3624) node {$v'_{12}$};
\draw [fill=black] (9.,4.) circle (3.0pt);
\draw[color=black] (9.1612,4.3644) node {$v'_6$};
\draw [fill=black] (18.,9.) circle (3.0pt);
\draw[color=black] (18.2142,9.4134) node {$v''_3$};
\draw [fill=black] (18.,-1.) circle (3.0pt);
\draw[color=black] (18.2142,-0.5746) node {$v''_1$};
\draw [fill=black] (19.,6.) circle (3pt);
\draw[color=black] (19.2042,6.4214) node {$v''_7$};
\draw [fill=black] (19.,-1.) circle (3.0pt);
\draw[color=black] (19.2042,-0.5746) node {$v''_4$};
\draw [fill=black] (21.,6.) circle (3.0pt);
\draw[color=black] (21.2062,6.4214) node {$v''_{11}$};
\draw [fill=black] (21.,1.) circle (3.0pt);
\draw[color=black] (21.2062,1.4274) node {$v''_8$};
\draw [fill=black] (19.,4.) circle (3.0pt);
\draw[color=black] (19.2042,4.4194) node {$v''_6$};
\draw [fill=black] (23.,4.) circle (3.0pt);
\draw[color=black] (23.3,4.4194) node {$v''_{12}$};
\draw [fill=black] (18.,2.) circle (3.0pt);
\draw[color=black] (18.2142,2.4174) node {$v''_2$};
\draw [fill=black] (21.,2.) circle (3.0pt);
\draw[color=black] (21.2062,2.4174) node {$v''_9$};
\draw [fill=black] (6.,4.) circle (3.0pt);
\draw[color=black] (6.25,4.3644) node {$v'_5$};
\draw [fill=black] (9.,3.) circle (3.0pt);
\draw[color=black] (9.3,3.4074) node {$v'_{10}$};
\draw [fill=black] (21.,4.) circle (3.0pt);
\draw[color=black] (21.3,4.4194) node {$v''_{10}$};
\draw [fill=black] (19.,2.) circle (3pt);
\draw[color=black] (19.2042,2.4174) node {$v''_5$};
\draw (6,-2) node[anchor=north west] {$\mathscr{T}_1$};
    		\draw (20,-2) node[anchor=north west] {$\mathscr{T}_2$};
\end{scriptsize}
\end{tikzpicture}
    \caption{\label{fig structurally isomorphic class} \textcolor{black}{Two T-meshes in $[\mathscr{T}]_{\mathrm I}$}}
\end{figure}

In the following, we provide the mathematical definition of dimensional stability of T-meshes for spline space $\mathbf{S}_{d}(\mathscr{T})$. The dimensional stability refers to the consistent dimension of the polynomial spline space across any two T-meshes within the same structurally isomorphic class.

\begin{definition}\label{def dimstable}
    \textcolor{black}{Let $\mathscr{T}$ be a T-mesh. We call $\mathscr{T}$ \textbf{dimensionally stable} for $\mathbf{S}_{d}$ if, for every $\mathscr{T}_1,\mathscr{T}_2\in[\mathscr{T}]_{\mathrm I}$, where $[\mathscr{T}]_{\mathrm I}\in\mathrm{SIC}$ is the structurally isomorphic class of $\mathscr{T}$, one has}
    \begin{equation}\label{eq:dimension equal}
        \dim \mathbf{S}_d(\mathscr{T}_1)=\dim \mathbf{S}_d(\mathscr{T}_2).
    \end{equation}
\end{definition}

The dimensional stability of a T-mesh $\mathscr{T}$ for a spline space $\mathbf{S}_{\mathbf{d}}^{\textcolor{black}{\mathbf{r}}}(\mathscr{T})$ depends on both the T-mesh $\mathscr{T}$ and its associated spline space $\mathbf{S}_{\mathbf{d}}^{\textcolor{black}{\mathbf{r}}}(\mathscr{T})$. As this paper focuses on the spline space $\mathbf{S}_{d}(\mathscr{T})$, the corresponding spline space is omitted when discussing the dimensional stability of $\mathscr{T}$ in subsequent sections for simplicity.

\begin{theorem}\label{thm dimension stable}
    \textcolor{black}{A T-mesh $\mathscr{T}$ is dimensionally stable if and only if, for any $\mathscr{T}_1, \mathscr{T}_2 \in [\mathscr{T}]_{\mathrm{I}}$,}
    $$\mathrm{rank}(M(T(\mathscr{T}_1))) = \mathrm{rank}(M(T(\mathscr{T}_2))).$$
\end{theorem}
\begin{proof}
    This conclusion is established from the dimension formula~\eqref{eqdim} immediately. Since the structurally isomorphic map preserves the numbers of cross-cuts and interior vertices, the dimension remains invariant if and only if the rank of the conformality matrix is constant.
    $\Box$
\end{proof}

In Theorem~\ref{thm dimension stable}, since $\mathscr{T}_1$ and $\mathscr{T}_2$ are chosen from the same structurally isomorphic class $[\mathscr{T}]_{\mathrm{I}}$, they satisfy the conditions in Definition~\ref{defsim}. The third condition of Definition~\ref{defsim} states that the structurally isomorphic map preserves the \textbf{mutual positions} of all horizontal and vertical $l$-edges. Because all vertices are formed by the intersections of these $l$-edges, the relative positions of the vertices are also determined by the configuration of the $l$-edges. Therefore, the rank condition $\mathrm{rank}(M(T(\mathscr{T}_1))) = \mathrm{rank}(M(T(\mathscr{T}_2)))$ depends intrinsically on the mutual positions of the vertices on $T(\mathscr{T})$.

In existing literature, such as \cite{Ins2011,Ins2012,tc2015,tc2017}, the dimensional stability is described as a property depending solely on the "topology" of a T-mesh without clarifying its precise meaning. By rigorously defining the dimensional stability, we find that the stability depends on both the connectivity of $l$-edges (the second condition of Definition~\ref{defsim}) and the mutual positions of vertices.

Furthermore, Theorem~\ref{thm dimension stable} reduces the global dimensional stability to the properties of the T-connected component $T(\mathscr{T})$, showing that stability depends on the connectivity of T $l$-edges and the mutual positions of vertices within $T(\mathscr{T})$. This motivates us to introduce a stronger type of stability in the next section, termed \textbf{absolute stability}, which depends purely on the connectivity of T $l$-edges and is independent of the mutual positions of vertices.

\subsection{Structurally similar map and absolute stability}\label{subsection: defintion ab}

As established by Theorem~\ref{thm dimension stable}, the dimensional stability of $\mathscr{T}$ involves both the connectivity of T $l$-edges and the mutual positions of vertices in $T(\mathscr{T})$. To further study the stability from a purely topological perspective, we introduce the concept of absolute stability. Unlike conventional dimensional stability, absolute stability requires the dimension of the spline space to remain constant under any variations that preserve only the connectivity of $T(\mathscr{T})$, regardless of the mutual positions of its vertices. Prior to discussing absolute stability, we first define the concept of \textbf{the generalized T-connected component.}

\begin{definition}[Generalized T-connected component]\label{def GT}
\color{black}
A \textbf{generalized T-connected component} is a finite connected edge--vertex configuration
\[
\mathrm{GT}=(L_h,L_v,V),
\]
with the following properties.
\begin{itemize}
    \item $L_h$ and $L_v$ are finite sets of closed horizontal and vertical segments, respectively. Distinct parallel segments are disjoint, and two perpendicular segments meet in at most one point. The elements of $L(\mathrm{GT})=L_h\cup L_v$ are called the \emph{$l$-edges of $\mathrm{GT}$}.
    \item Each $l\in L(\mathrm{GT})$ carries a finite set $V(l)\subset l$ of marked vertices containing both endpoints of $l$. Whenever $l_i\cap l_j\ne\varnothing$, the intersection point belongs to both $V(l_i)$ and $V(l_j)$, and $V=\bigcup_{l\in L(\mathrm{GT})}V(l)$.
    \item Connectedness means that any two $l$-edges can be joined by a finite chain of pairwise intersecting $l$-edges.
\end{itemize}
A vertex belonging to two $l$-edges is a multi-vertex; a vertex belonging to exactly one $l$-edge is a mono-vertex. We write $n(l)=|V(l)|$ and $m(l)$ for the number of mono-vertices on $l$.
\end{definition}

{\color{black}A generalized T-connected component is not required to be associated with, or realizable as the T-connected component of, any T-mesh. Every actual $T(\mathscr{T})$ is a generalized T-connected component, but the converse need not hold. The non-intersection vertices may move along their assigned $l$-edges, provided that they remain distinct and the prescribed incidence relations are retained. The conformality conditions, conformality matrix, conformality vector space, partitions, and relative mono-vertices are defined for $\mathrm{GT}$ exactly as for $T(\mathscr{T})$, using the marked vertices and their coordinates.}

 \begin{example}\label{exm gTc}
        \color{black}
        Consider the T-mesh $\mathscr{T}$ shown on the right of Fig.~\ref{fig ssi}. The left panel of Fig.~\ref{fig difference between GT and T} is its actual T-connected component $T(\mathscr{T})$. The right panel is a separate generalized T-connected component illustrating Definition~\ref{def GT}; it is not asserted to be associated with $\mathscr{T}$ or with any other T-mesh.
    \end{example}
    \begin{figure}
        \centering
        \begin{tikzpicture}[line cap=round,line join=round,>=triangle 45,x=1.0cm,y=1.0cm]
\draw [line width=1.pt] (1.,3.)-- (5.,3.);
\draw [line width=1.pt] (2.,3.)-- (2.,6.);
\draw [line width=1.pt] (0.,4.)-- (3.,4.);
\draw [line width=1.pt] (3.,5.)-- (3.,1.);
\draw [line width=1.pt] (7.,3.)-- (12.,3.);
\draw [line width=1.pt] (8.,4.)-- (11.,4.);
\draw [line width=1.pt] (9.,5.)-- (9.,1.);
\draw [line width=1.pt] (10.,6.)-- (10.,2.);
\draw (2.08,0.92) node[anchor=north west] {$T(\mathscr{T})$};
\draw (9.18,0.84) node[anchor=north west] {$\mathrm{GT}$};
\begin{scriptsize}
\draw [fill=black] (1.,3.) circle (2.0pt);
\draw[color=black] (1.14,3.25) node {$v_8$};
\draw [fill=black] (5.,3.) circle (2.0pt);
\draw[color=black] (5.14,3.25) node {$v_{11}$};
\draw [fill=black] (2.,3.) circle (2.0pt);
\draw[color=black] (2.14,3.25) node {$v_9$};
\draw [fill=black] (2.,6.) circle (2.0pt);
\draw[color=black] (2.14,6.25) node {$v_1$};
\draw [fill=black] (0.,4.) circle (2.0pt);
\draw[color=black] (0.14,4.25) node {$v_4$};
\draw [fill=black] (3.,4.) circle (2.0pt);
\draw[color=black] (3.14,4.25) node {$v_7$};
\draw [fill=black] (3.,5.) circle (2.0pt);
\draw[color=black] (3.14,5.25) node {$v_3$};
\draw [fill=black] (3.,1.) circle (2.0pt);
\draw[color=black] (3.25,1.25) node {$v_{12}$};
\draw [fill=black] (1.,4.) circle (2.0pt);
\draw[color=black] (1.14,4.25) node {$v_5$};
\draw [fill=black] (2.,4.) circle (2.0pt);
\draw[color=black] (2.14,4.25) node {$v_6$};
\draw [fill=black] (2.,5.) circle (2.0pt);
\draw[color=black] (2.14,5.25) node {$v_2$};
\draw [fill=black] (3.,3.) circle (2.0pt);
\draw[color=black] (3.25,3.25) node {$v_{10}$};
\draw [fill=black] (7.,3.) circle (2.0pt);
\draw[color=black] (7.14,3.25) node {$v'_8$};
\draw [fill=black] (12.,3.) circle (2.0pt);
\draw[color=black] (12.14,3.25) node {$v'_{11}$};
\draw [fill=black] (8.,4.) circle (2.0pt);
\draw[color=black] (8.14,4.25) node {$v'_4$};
\draw [fill=black] (11.,4.) circle (2.0pt);
\draw[color=black] (11.14,4.25) node {$v'_5$};
\draw [fill=black] (9.,5.) circle (2.0pt);
\draw[color=black] (9.14,5.25) node {$v'_1$};
\draw [fill=black] (9.,1.) circle (2.0pt);
\draw[color=black] (9.14,1.25) node {$v'_2$};
\draw [fill=black] (10.,6.) circle (2.0pt);
\draw[color=black] (10.14,6.25) node {$v'_3$};
\draw [fill=black] (10.,2.) circle (2.0pt);
\draw[color=black] (10.25,2.25) node {$v'_{12}$};
\draw [fill=black] (9.,4.) circle (2.0pt);
\draw[color=black] (9.14,4.33) node {$v'_6$};
\draw [fill=black] (10.,4.) circle (2.0pt);
\draw[color=black] (10.14,4.33) node {$v'_7$};
\draw [fill=black] (9.,3.) circle (2.0pt);
\draw[color=black] (9.14,3.33) node {$v'_9$};
\draw [fill=black] (10.,3.) circle (2.0pt);
\draw[color=black] (10.25,3.33) node {$v'_{10}$};
\end{scriptsize}
\end{tikzpicture}
        \caption{\label{fig difference between GT and T}The difference between $T(\mathscr{T})$ and $\mathrm{GT}$.}
    \end{figure}

    We next define the concept of a \textbf{structurally similar map}, which preserves the connectivity of edges within the generalized T-connected component.

\begin{definition}[Structurally similar map]\label{def ssm}
\color{black}
Let $\mathrm{GT}_1$ and $\mathrm{GT}_2$ be generalized T-connected components with $l$-edge sets $L(\mathrm{GT}_1)$ and $L(\mathrm{GT}_2)$. A \textbf{structurally similar map} is a bijection
\[
\beta:L(\mathrm{GT}_1)\longrightarrow L(\mathrm{GT}_2)
\]
such that, for every $l\in L(\mathrm{GT}_1)$ and every pair $l_i,l_j\in L(\mathrm{GT}_1)$,
\[
n(l)=n(\beta(l))
\quad\text{and}\quad
l_i\cap l_j\ne\varnothing
\Longleftrightarrow
\beta(l_i)\cap\beta(l_j)\ne\varnothing.
\]
No condition is imposed on the horizontal/vertical orientation or on the relative coordinates of corresponding $l$-edges.
\end{definition}

\begin{example}\label{exm:gTc}
    \color{black}
    The left panel of Fig.~\ref{fig difference between GT} shows the particular generalized T-connected component $\mathrm{GT}_1=T(\mathscr{T}_1)$ obtained from the T-mesh $\mathscr{T}_1$ in Fig.~\ref{fig:tmesh_example}. The right panel shows one particular generalized T-connected component $\mathrm{GT}_2$. A T-mesh determines the specific member $T(\mathscr{T})$, while its structurally similar class generally contains many other generalized realizations; therefore, no uniqueness of a generalized component ``corresponding to'' a T-mesh is intended.

    The displayed $\mathrm{GT}_1$ and $\mathrm{GT}_2$ are not structurally similar. Indeed, the $l$-edges $v_1v_2$ and $v_7v_8$ do not intersect in $\mathrm{GT}_1$, whereas their proposed counterparts $v'_1v'_2$ and $v'_7v'_8$ intersect in $\mathrm{GT}_2$. Hence the incidence condition in Definition~\ref{def ssm} fails.
\end{example}

\begin{figure}[htbp]
    \centering
    \begin{tikzpicture}[scale=0.6, line cap=round, line join=round, every node/.style={font=\small}]
        \draw [line width=1pt, black] (1,1) rectangle (9,9);
        \draw [line width=1pt, black] (1,8) -- (9,8);
        \draw [line width=1pt, black] (1,7) -- (9,7);
        \draw [line width=1pt, black] (1,2) -- (9,2);
        \draw [line width=1pt, black] (2,1) -- (2,9);
        \draw [line width=1pt, black] (8,1) -- (8,9);

        \draw [line width=1pt, blue] (2,6) -- (6,6); 
        \draw [line width=1pt, black] (4,6) -- (4,9);
        \draw [line width=1pt, blue] (3,2) -- (3,7); 
        \draw [line width=1pt, blue] (7,2) -- (7,8); 
        \draw [line width=1pt, blue] (2,4) -- (7,4); 
        \draw [line width=1pt, black] (6,4) -- (6,9);

        \filldraw [black] (2,6) circle (2.5pt) node[above left] {$v_1$};
        \filldraw [black] (6,6) circle (2.5pt) node[above right] {$v_2$};
        \filldraw [black] (3,7) circle (2.5pt) node[above left] {$v_3$};
        \filldraw [black] (3,2) circle (2.5pt) node[below left] {$v_4$};
        \filldraw [black] (2,4) circle (2.5pt) node[above left] {$v_5$};
        \filldraw [black] (7,4) circle (2.5pt) node[below right] {$v_6$};
        \filldraw [black] (7,8) circle (2.5pt) node[above right] {$v_7$};
        \filldraw [black] (7,2) circle (2.5pt) node[below right] {$v_8$};
    \end{tikzpicture}

    \caption{\label{fig:tmesh_example}A T-mesh $\mathscr{T}_1$.}
\end{figure}

\begin{figure}[htbp]
    \centering
\begin{tikzpicture}[scale=0.75, line cap=round, line join=round, >=triangle 45, every node/.style={font=\small}]
    \draw [line width=1pt, blue] (2,6) -- (6,6); 
    \draw [line width=1pt, blue] (3,2) -- (3,7); 
    \draw [line width=1pt, blue] (2,4) -- (7,4); 
    \draw [line width=1pt, blue] (7,2) -- (7,8); 

    \filldraw [black] (2,6) circle (2.0pt) node[above left] {$v_1$};
    \filldraw [black] (6,6) circle (2.0pt) node[above right] {$v_2$};
    \filldraw [black] (3,7) circle (2.0pt) node[above] {$v_3$};
    \filldraw [black] (3,2) circle (2.0pt) node[below] {$v_4$};
    \filldraw [black] (2,4) circle (2.0pt) node[left] {$v_5$};
    \filldraw [black] (7,4) circle (2.0pt) node[below right] {$v_6$};
    \filldraw [black] (7,8) circle (2.0pt) node[above] {$v_7$};
    \filldraw [black] (7,2) circle (2.0pt) node[below] {$v_8$};

    \filldraw [black] (3,6) circle (2.0pt); 
    \filldraw [black] (4,6) circle (2.0pt); 
    \filldraw [black] (3,4) circle (2.0pt); 
    \filldraw [black] (6,4) circle (2.0pt); 
    \filldraw [black] (7,7) circle (2.0pt); 

    \node at (4.5, 0.3) {(a) $\mathrm{GT}_1$};

    \draw [line width=1pt, black] (9,4) -- (14,4); 
    \draw [line width=1pt, black] (10,6) -- (13,6); 
    \draw [line width=1pt, black] (11,2) -- (11,7); 
    \draw [line width=1pt, black] (12,3) -- (12,8); 

    \filldraw [black] (9,4) circle (2.0pt) node[left] {$v'_5$};
    \filldraw [black] (14,4) circle (2.0pt) node[right] {$v'_6$};
    \filldraw [black] (11,7) circle (2.0pt) node[above] {$v'_3$};
    \filldraw [black] (11,2) circle (2.0pt) node[below] {$v'_4$};

    \filldraw [black] (10,6) circle (2.0pt) node[left] {$v'_1$};
    \filldraw [black] (13,6) circle (2.0pt) node[right] {$v'_2$};
    \filldraw [black] (12,8) circle (2.0pt) node[above] {$v'_7$};
    \filldraw [black] (12,3) circle (2.0pt) node[below] {$v'_8$};

    \filldraw [black] (11,4) circle (2.0pt);
    \filldraw [black] (12,4) circle (2.0pt);
    \filldraw [black] (11,6) circle (2.0pt);
    \filldraw [black] (12,6) circle (2.0pt);

    \node at (11.5, 0.3) {(b) $\mathrm{GT}_2$};
\end{tikzpicture}
    \caption{\label{fig difference between GT}Two generalized T-connected components $\mathrm{GT}_1$ and $\mathrm{GT}_2$ that are not structurally similar}
\end{figure}

\textcolor{black}{The proof of Proposition~\ref{proof of equivalence of sim} applies verbatim to give the following result.}

        \begin{proposition}
           \textcolor{black}{Let $\mathcal{G}$ be the set of all generalized T-connected components. Structural similarity is an equivalence relation on $\mathcal{G}$.}
        \end{proposition}

        The concept of \textbf{structurally similar class} is defined similarly as follows.

\begin{definition}\label{def ssc}
    \color{black}
    Let $\mathcal{G}$ denote the set of all generalized T-connected components. An equivalence class induced by structural similarity is called a \textbf{structurally similar class}. The set of all such classes is denoted by $\mathrm{SSC}(\mathcal{G})$, or simply $\mathrm{SSC}$.
\end{definition}

{\color{black}For a fixed component $T(\mathscr{T})$ with $l$-edges $l_1,\ldots,l_t$, a member of $[T(\mathscr{T})]_{\mathrm S}$ is constructed by choosing $t$ axis-aligned segments $\widetilde l_1,\ldots,\widetilde l_t$, placing exactly $n(l_i)$ marked vertices on $\widetilde l_i$, and imposing
\[
\widetilde l_i\cap\widetilde l_j\ne\varnothing
\quad\Longleftrightarrow\quad
l_i\cap l_j\ne\varnothing
\]
for every pair. Non-intersection vertices may otherwise be placed arbitrarily on their assigned segments. The class is generally infinite; the quantifier over all members of $[T(\mathscr{T})]_{\mathrm S}$ is therefore a mathematical quantifier over all such realizations, not a request for a finite enumeration.}

All generalized T-connected components within one class have the same incidence structure and the same number of vertices on corresponding $l$-edges. Based on $\mathrm{SSC}$, we introduce the stronger notion of dimensional absolute stability.

    \begin{definition}\label{def dimensionally absolutely stable}
        \textcolor{black}{Let $\mathscr{T}$ be a T-mesh. We call $\mathscr{T}$ \textbf{dimensionally absolutely stable} for $\mathbf{S}_{d}$ if, for every $\mathrm{GT}_1,\mathrm{GT}_2\in[T(\mathscr{T})]_{\mathrm S}$, where $[T(\mathscr{T})]_{\mathrm S}\in\mathrm{SSC}$ is the structurally similar class of $T(\mathscr{T})$, one has}
    $$\dim\mathsf{CVS}[\mathrm{GT}_1]=\dim\mathsf{CVS}[\mathrm{GT}_2],$$
    where $\mathsf{CVS}[\mathrm{GT}_i]$ is the conformality vector space of $\mathrm{GT}_i$ for $i=1,2$.
    \end{definition}

   Dimensional absolute stability depends solely on the topology of $T(\mathscr{T})$. By the definition of dimensional absolute stability, the \textcolor{black}{following} theorem is immediate.

    \begin{theorem}\label{thm as to s}
      If ${\mathscr{T}}$ is a dimensionally absolutely stable T-mesh, then $\mathscr{T}$ is a dimensionally stable T-mesh.
    \end{theorem}
    \begin{proof}
\color{black}
Let $\mathscr{T}_1,\mathscr{T}_2\in[\mathscr{T}]_{\mathrm I}$ and let $\alpha:L(\mathscr{T}_1)\to L(\mathscr{T}_2)$ be a structurally isomorphic map. Its restriction
\[
\beta=\alpha|_{L(T(\mathscr{T}_1))}:L(T(\mathscr{T}_1))\longrightarrow L(T(\mathscr{T}_2))
\]
is a bijection. The type-preserving condition maps T $l$-edges to T $l$-edges, and incidence preservation gives
\[
l_i\cap l_j\ne\varnothing
\Longleftrightarrow
\beta(l_i)\cap\beta(l_j)\ne\varnothing.
\]
The structural isomorphism also preserves the vertices lying on corresponding $l$-edges, so $n(l)=n(\beta(l))$. Hence $\beta$ is structurally similar, and both $T(\mathscr{T}_1)$ and $T(\mathscr{T}_2)$ lie in $[T(\mathscr{T})]_{\mathrm S}$.

Absolute stability therefore yields
\[
\dim\mathsf{CVS}[T(\mathscr{T}_1)]
=
\dim\mathsf{CVS}[T(\mathscr{T}_2)],
\]
or equivalently, equality of the ranks of their conformality matrices. Theorem~\ref{thm dimension stable} then implies that $\mathscr{T}$ is dimensionally stable.

$\Box$
\end{proof}

    We conclude with a remark on the generalized T-connected component. A generalized T-connected component is not necessarily a T-connected component of a T-mesh. For example, in Fig.~\ref{fig difference between GT and T}, $\mathrm{GT}$ cannot be embedded in any T-mesh, confirming it is not a T-connected component of any T-mesh. However, in certain cases, such as those shown in \cite{Ins2011,tc2015}, it may be.


    \begin{example}\label{exm ssi}
    \color{black}
    Fig.~\ref{fig ssi} presents two T-meshes: $\mathscr{T}_1$ from \cite{tc2015} and $\mathscr{T}_2$ from \cite{Ins2011}. Their T-connected components $T(\mathscr{T}_1)$ and $T(\mathscr{T}_2)$ are structurally similar, although the two full T-meshes differ in features outside those components, including their cross-cuts and rays.
    \end{example}
    \begin{figure}
    	\centering
    	\begin{tikzpicture}[line cap=round,line join=round,>=triangle 45,x=1.0cm,y=1.0cm,scale=1.35]
    		\draw [line width=1.pt] (1.,1.)-- (1.,6.);
    		\draw [line width=1.pt] (1.,1.5)-- (6,1.5);
    		\draw [line width=1.pt] (1.,1.)-- (6.,1.);
    		\draw [line width=1.pt] (6.,1.)-- (6.,6.);
    		\draw [line width=1.pt] (1.,6.)-- (6.,6.);
    		\draw [line width=1.pt] (1.,5.5)-- (6.,5.5);
    		\draw [line width=1.pt] (1.4,6.)-- (1.4,1.);
    		\draw [line width=1.pt] (5.6,6.)-- (5.6,1.);
    		\draw [line width=1.pt] (3,5.5)-- (3.,3.);
    		\draw [line width=1.pt] (1.4,3)-- (4.,3.);
    		\draw [line width=1.pt] (3,4)-- (5.6,4);
    		\draw [line width=1.pt] (4,4)-- (4,1.5);
    		\draw [line width=1.pt] (1.,3.4)-- (3,3.4);
    		\draw [line width=1.pt] (3.6,6.)-- (3.6,4);
    		\draw [line width=1.pt] (3.4,3)-- (3.4,1.);
    		\draw [line width=1.pt] (4,3.6)-- (6.,3.6);
    		\draw [line width=1.pt] (8.,6.)-- (8.,1.);
    		\draw [line width=1.pt] (8.,1.)-- (13.,1.);
    		\draw [line width=1.pt] (8.,6.)-- (13.,6.);
    		\draw [line width=1.pt] (13.,6.)-- (13.,1.);
    		\draw [line width=1.pt] (9.,6.)-- (9.,1.);
    		\draw [line width=1.pt] (8.,5.)-- (13.,5.);
    		\draw [line width=1.pt] (8.,5.6)-- (13.,5.6);
    		\draw [line width=1.pt] (8.5,6.)-- (8.5,1.);
    		\draw [line width=1.pt] (8.,1.54)-- (13.,1.54);
    		\draw [line width=1.pt] (12.4,6.)-- (12.4,1.);
    		\draw [line width=1.pt] (10,5.6)-- (10.,3.);
    		\draw [line width=1.pt] (9.,3.)-- (12.4,3);
    		\draw [line width=1.pt] (11.,5.)-- (11.,1.5);
    		\draw [line width=1.pt] (8.5,4)-- (11.,4.);
    		\draw (3.2485671561654423,0.7144552129249617) node[anchor=north west] {$\mathscr{T}_1$};
    		\draw (10.426558351570758,0.7144552129249617) node[anchor=north west] {$\mathscr{T}_2$};
    		\begin{scriptsize}
    			\draw [fill=black] (3,5.5) circle (2.0pt);
    			\draw[color=black] (3.12,5.7) node {$v_1$};
    			\draw [fill=black] (3.,3.) circle (2.0pt);
    			\draw[color=black] (3.12,3.2) node {$v_2$};
    			\draw [fill=black] (1.4,3) circle (2.0pt);
    			\draw[color=black] (1.52,3.2) node {$v_3$};
    			\draw [fill=black] (3.4,3.) circle (2.0pt);
    			\draw[color=black] (3.52,3.2) node {$v_4$};

                \draw [fill=black] (5.6,4) circle (2.0pt);
    			\draw[color=black] (5.72,4.2) node {$v_8$};

    			\draw [fill=black] (4,3.6) circle (2.0pt);
    			\draw[color=black] (4.15,3.7) node {$v_6$};
    			\draw [fill=black] (4,4) circle (2.0pt);
    			\draw[color=black] (4.12,4.2) node {$v_7$};
    			\draw [fill=black] (4,1.5) circle (2.0pt);
    			\draw[color=black] (4.12,1.7) node {$v_5$};
    			\draw [fill=black] (10,5.56) circle (2.0pt);
    			\draw[color=black] (10.136,5.75) node {$v_1'$};
    			\draw [fill=black] (10.,3.) circle (2.0pt);
    			\draw[color=black] (10.136,3.2) node {$v_2'$};
    			\draw [fill=black] (9.,3.) circle (2.0pt);
    			\draw[color=black] (9.136,3.2) node {$v_3'$};
    			\draw [fill=black] (12.398414031508805,2.994873892813983) circle (2.0pt);
    			\draw[color=black] (12.55,3.2) node {$v_4'$};
    			\draw [fill=black] (11.,5.) circle (2.0pt);
    			\draw[color=black] (11.141525916399099,5.15) node {$v_6'$};
    			\draw [fill=black] (11.,1.54807) circle (2.0pt);
    			\draw[color=black] (11.141525916399099,1.7224) node {$v_5'$};
    			\draw [fill=black] (8.504207257494103,4) circle (2.0pt);
    			\draw[color=black] (8.635599996109669,4.2) node {$v_8'$};
    			\draw [fill=black] (11.,4.) circle (2.0pt);
    			\draw[color=black] (11.141525916399099,4.2) node {$v_7'$};
    		\end{scriptsize}
    	\end{tikzpicture}
    	\caption{\label{fig ssi} \textcolor{black}{Two T-meshes with structurally similar T-connected components}}
    \end{figure}

\subsection{Summary of Section 3}
 This section introduces the structurally isomorphic map between two T-meshes, which preserves the connectivity and mutual positions of $l$-edges. Based on the structurally isomorphic class, the dimensional stability of a T-mesh is defined as an invariant property of the spline space dimension. We show that dimensional stability is equivalently determined by the rank stability of the conformality matrix associated with the T-connected component $T(\mathscr{T})$. Because this conventional dimensional stability intrinsically depends on both the topological connectivity of $l$-edges and the mutual positions of vertices within $T(\mathscr{T})$, we further introduce a stronger notion, termed absolute stability, via structurally similar maps. Absolute stability requires the dimension to depend purely on the mesh connectivity and to be completely independent of the relative position configuration of vertices.

Fig.~\ref{fig:stability_map_relations} depicts the relationships among these concepts. For a T-mesh $\mathscr{T}$, absolute stability (defined over the structurally similar class) implies conventional dimensional stability (defined over the structurally isomorphic class). The relationships among these classes are described as follows:
 \[\textcolor{black}{\mathscr{T}_1,\mathscr{T}_2\in[\mathscr{T}]_{\mathrm I}\quad\Longrightarrow\quad T(\mathscr{T}_1),T(\mathscr{T}_2)\in[T(\mathscr{T})]_{\mathrm S}.}\]

\begin{figure}[ht]
\centering
\begin{tikzcd}
\text{dimensionally absolutely stable} \arrow[d, Rightarrow] \arrow[r, Leftrightarrow]
    & \text{Structurally similar class} \\
\text{dimensionally stable} \arrow[r, Leftrightarrow]
    & \text{Structurally isomorphic class} \arrow[u, Rightarrow]
\end{tikzcd}
\caption{\textcolor{black}{Relationship between dimensional stability types and structural classes.} Vertical arrows indicate that a stronger stability condition or class implies a weaker one, while horizontal arrows denote equivalence.}
\label{fig:stability_map_relations}
\end{figure}

\textcolor{black}{The concepts presented in Section~3 yield the following preliminary classification of T-meshes by dimensional stability.}
\begin{equation}
\text{T-meshes} \begin{dcases*}
    \text{dimensionally unstable T-meshes} \\
    \text{dimensionally stable T-meshes} \begin{dcases*}
        \text{dimensionally absolutely stable T-meshes}\\
        \text{dimensionally non-absolutely stable T-meshes}
    \end{dcases*}
\end{dcases*}
\label{eq:Tmesh-classification}
\end{equation}

Finally, we present the new symbols that appeared in this section in Table~\ref{tab:symbols section 3}.
\begin{table}[H]
\centering
\caption{List of Symbols in Section~\ref{sec:Structurally Similar Map}}
\begin{tabular}{|c|l|}
\hline
\textbf{Symbol} & \textbf{Explanation} \\
\hline
$L(\mathscr{T})$ & The set of all $l$-edges (including boundary $l$-edges) of $\mathscr{T}$ \\
\hline
$\alpha$ & The structurally isomorphic map \\
\hline
$\mathrm{cor}(l)$ & The coordinate value of $l$-edge $l$\\
\hline
$L_1(\mathscr{T})$ & The set of all horizontal $l$-edges \\
\hline
$L_2(\mathscr{T})$ & The set of all vertical $l$-edges \\
\hline
$U$ & The set of all T-meshes\\
\hline
$\sim_{\mathrm{I}}$ & The structurally isomorphic equivalence relation in $U$\\
\hline
$\mathrm{SIC}$ & The set of all structurally isomorphic classes\\
\hline
$\mathrm{GT}$ & The generalized T-connected component.\\
\hline
$\beta$ & The structurally similar map.\\
\hline
$\mathrm{SSC}$ & The set of all structurally similar classes\\
\hline
\end{tabular}
\label{tab:symbols section 3}
\end{table}

\section{The Existence and Uniqueness of the Complete Partition}\label{sec: complete partition}

{\color{black}The previous section defined dimensional stability. This section establishes the theoretical basis of the complete partition introduced in \cite{huang2023}. The partition separates a T-connected component into a diagonalizable part and a completely non-diagonalizable component (CNDC), thereby reducing dimensional-stability analysis to the rank stability of the conformality matrix associated with the CNDC.}

The primary focus of this section is the theoretical foundation of the CNDC, namely its existence and uniqueness. This theoretical necessity arises because Algorithm~\ref{alg} provides a decomposition by sequentially extracting diagonalizable T $l$-edges. A critical question is whether different choices of diagonalizable T $l$-edges during the extraction process lead to different final CNDCs. In other words, the formulation of Algorithm~\ref{alg} implies a potential dependence on the extraction order. This section thoroughly resolves this issue by proving that the complete partition is entirely independent of the extraction order, meaning the resulting CNDC is unique.

To achieve this, we first introduce the $k$-partition, a generalization of the bipartite partition and $t$-partitions. We then establish a key property of the $k$-partition, which directly leads to the proof of existence and uniqueness. Finally, based on this uniqueness, we derive a series of results regarding the dimension formula.

\subsection{The definition of \texorpdfstring{$k$}{k}-partition}

\begin{definition}[$k$-partition]\label{def4.1}
\color{black}
Let $T(\mathscr{T})$ have T $l$-edges $l_1,\ldots,l_t$. Fix an ordering
\[
l_1\succ l_2\succ\cdots\succ l_t
\]
and indices
\[
0=m_0<m_1<\cdots<m_k=t.
\]
For $i=1,\ldots,k$, let
\[
\mathcal{L}_i=\{l_{m_{i-1}+1},\ldots,l_{m_i}\}.
\]
The ordering and the indices determine an indexed collection $\{\overline T_1,\ldots,\overline T_k\}$ as follows: $\overline T_i$ contains precisely the $l$-edges in $\mathcal{L}_i$, and every vertex shared with an $l$-edge in an earlier block $\mathcal{L}_1\cup\cdots\cup\mathcal{L}_{i-1}$ is omitted from $\overline T_i$. A collection obtained from some ordering and some indices in this manner is called \emph{a $k$-partition} of $T(\mathscr{T})$.
\end{definition}

{\color{black}The ordering $l_1\succ\cdots\succ l_t$ and the break indices $m_0,\ldots,m_k$ are part of the data of a $k$-partition. Thus a given T-connected component generally admits more than one $k$-partition. Each T $l$-edge belongs to exactly one component, and a vertex shared by different blocks is assigned to the earliest block containing one of its incident $l$-edges.}

\begin{figure}[ht]
	\centering
	\vspace{-0.35cm}
	\subfigtopskip=2pt
	\subfigbottomskip=2pt
	\subfigcapskip=-5pt
	\subfigure[A T-mesh and its T-connected component]{
		\label{fig2.sub.1}
		\includegraphics[width=0.4\linewidth]{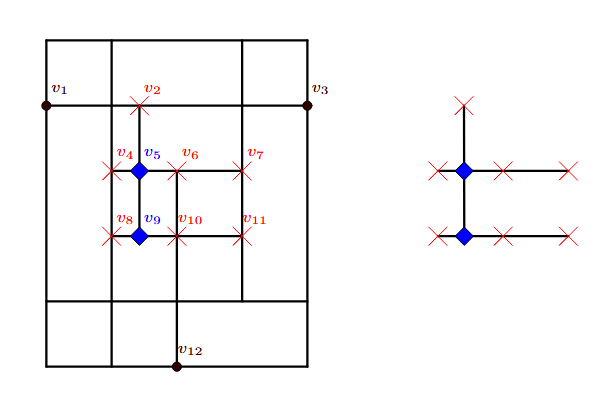}}
	\quad
	\subfigure[A $2$-partition of T-connected component]{
		\label{fig2.sub.2}
		\includegraphics[width=0.5\linewidth]{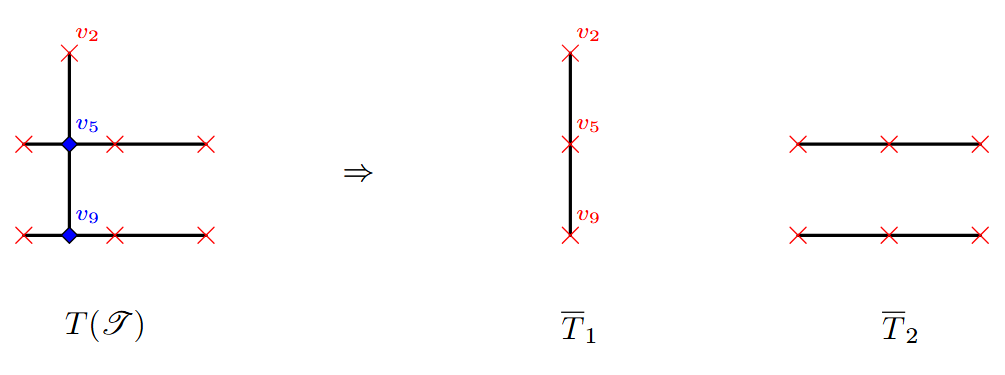}}
	\caption{\label{fig2}\textcolor{black}{An example of a $2$-partition of a T-connected component}}
\end{figure}

\begin{example}\label{exm3.1}
    Consider the T-mesh $\mathscr{T}$ in Fig.~\ref{fig2.sub.1}. The T-connected component contains three T $l$-edges: $v_2v_9$, $v_4v_7$, and $v_8v_{11}$.

    \textcolor{black}{Choose the ordering $v_2v_9\succ v_4v_7\succ v_8v_{11}$ and the break index $m_1=1$. Then $\mathcal{L}_1=\{v_2v_9\}$ and $\mathcal{L}_2=\{v_4v_7,v_8v_{11}\}$, yielding the $2$-partition $\{\overline{T}_1, \overline{T}_2\}$ shown in Fig.~\ref{fig2.sub.2}.} The intersection vertices $v_5$ and $v_9$ are fully assigned to $\overline{T}_1$. Consequently, they act as mono-vertices within $\overline{T}_1$, while the remaining edges in $\overline{T}_2$ are stripped of these interior vertices.

    {\color{black}Conversely, if the block order is reversed, the component containing the two horizontal T $l$-edges retains $v_5$ and $v_9$, while the later component containing $v_2v_9$ omits these shared vertices. Thus different choices of the ordering can produce different $2$-partitions, even though the underlying T-connected component is unchanged.}
\end{example}

\begin{definition}[Relative mono-vertex]\label{def4.2}
\color{black}
Let $S$ be a component of a $k$-partition and let $l$ be an $l$-edge retained in $S$. The notation $l|_S$ means the $l$-edge $l$ together with the vertices retained on it in $S$. Define
\[
V_{\mathrm m}(l|_S)=\{v\in V(l|_S): v\text{ belongs to no other $l$-edge of }S\}.
\]
The vertices in $V_{\mathrm m}(l|_S)$ are called the \textbf{relative mono-vertices of $l$ with respect to $S$}, and
\[
m(l|_S)=|V_{\mathrm m}(l|_S)|.
\]
\end{definition}

{\color{black}For example, in Fig.~\ref{fig2.sub.2}, the vertices $v_5$ and $v_9$ are multi-vertices in the full component $T(\mathscr{T})$, but they are relative mono-vertices of $v_2v_9$ in the isolated component $\overline T_1$. We use $S\subseteq S'$ to mean that every $l$-edge and every retained vertex of the edge--vertex structure $S$ also belongs to $S'$.}

\begin{lemma}\label{lem4.1}
\color{black}
Let $\{\overline T_1,\ldots,\overline T_k\}$ and $\{\overline T'_1,\ldots,\overline T'_k\}$ be two arbitrary $k$-partitions of $T(\mathscr{T})$. If $\overline T_i\subseteq\overline T'_j$, then, for every T $l$-edge $l\in\overline T_i$,
\[
m(l|_{\overline T'_j})\le m(l|_{\overline T_i}).
\]
\end{lemma}
\begin{proof}
\color{black}
Let $v\in V_{\mathrm m}(l|_{\overline T'_j})$. By definition, no $l$-edge of $\overline T'_j$ other than $l$ contains $v$. Since $\overline T_i\subseteq\overline T'_j$, no $l$-edge of $\overline T_i$ other than $l$ contains $v$ either. Hence
\[
V_{\mathrm m}(l|_{\overline T'_j})\subseteq V_{\mathrm m}(l|_{\overline T_i}).
\]
Taking cardinalities gives the stated inequality.

$\Box$
\end{proof}

Lemma~\ref{lem4.1} establishes an intrinsic geometric monotonicity under sub-structure restriction: as a partition component expands, it incorporates more surrounding edges, which potentially reintroduces intersection constraints and thereby monotonically reduces the number of relative mono-vertices. This property provides the essential theoretical bridge to establish the order-independence of Algorithm~\ref{alg} and the uniqueness of the complete partition.

\begin{example}
    Consider the 2-partition in Example~\ref{exm3.1} (Fig.~\ref{fig2.sub.2}). Let the sub-component $\overline{T}'_1 = T(\mathscr{T})$ and $\overline{T}_1 = \{v_2v_9\}$, which clearly satisfies $\overline{T}_1 \subseteq \overline{T}'_1$. For the edge $v_2v_9$, it has only 1 relative mono-vertex in the global mesh $T(\mathscr{T})$ (which is $v_2$). However, it has 3 relative mono-vertices in $\overline{T}_1$ ($v_2, v_5, v_9$). Thus, $m(v_2v_9|_{\overline{T}'_1}) = 1 < 3 = m(v_2v_9|_{\overline{T}_1})$, directly verifying Lemma~\ref{lem4.1}.
\end{example}

\subsection{The existence and uniqueness of CNDC}

In this subsection, we prove the existence and uniqueness of the CNDC. This proof is necessary for two reasons. First, while Algorithm~\ref{alg} offers a practical procedure for decomposing a T-connected component, it relies on sequentially extracting diagonalizable T $l$-edges. It remains unverified whether different extraction orders could yield structurally distinct final CNDCs; thus, proving order-independence is crucial to ensure the algorithm's validity. Second, as established in \cite{huang2023}, the computation of the spline dimension formula ultimately reduces to analyzing the rank of the conformality matrix associated with the CNDC. This dimension reduction framework strictly hinges on the prerequisite that the CNDC is uniquely determined by the intrinsic topology of the T-mesh. To address these issues, we present the main theorem below, which is proved by applying Lemma~\ref{lem4.1}.

Next, we present the theorem on the existence and uniqueness of the complete partition.

\begin{theorem}[Existence and Uniqueness of the Complete Partition]\label{thmcompuni}
\color{black}
Let $T(\mathscr{T})$ be the T-connected component of a T-mesh $\mathscr{T}$. Then there exists a unique complete partition
\[
\{\overline T_{\mathrm{CNDC}},\overline T_{\mathrm D}\}
\]
of $T(\mathscr{T})$.
\end{theorem}

\begin{proof}
\color{black}
Set $\Omega_0=T(\mathscr{T})$. Recursively define
\[
D_i=\{l\in\Omega_{i-1}:m(l|_{\Omega_{i-1}})\ge d+1\},
\qquad
\Omega_i=\Omega_{i-1}\setminus D_i.
\]
Because $T(\mathscr{T})$ has finitely many T $l$-edges, there is a first index $q$ for which $D_{q+1}=\varnothing$. Thus every $l\in\Omega_q$ satisfies $m(l|_{\Omega_q})<d+1$, while the edges in $D_1\cup\cdots\cup D_q$ can be removed layer by layer in an order compatible with the displayed construction. This proves existence, with
\[
\overline T_{\mathrm{CNDC}}=\Omega_q,
\qquad
\overline T_{\mathrm D}=\bigcup_{i=1}^{q}D_i.
\]

It remains to prove that an arbitrary sequential execution of Algorithm~\ref{alg} has the same output. Let $C$ be its terminal component and let $D$ be the set of extracted T $l$-edges.

First, every intrinsic layer $D_i$ is contained in $D$. We argue by induction. Suppose that $D_1,\ldots,D_{i-1}\subseteq D$. Then $C\subseteq\Omega_{i-1}$. If some $l\in D_i$ remained in $C$, Lemma~\ref{lem4.1} would give
\[
m(l|_C)\ge m(l|_{\Omega_{i-1}})\ge d+1,
\]
contradicting termination of the algorithm. Hence $D_i\subseteq D$.

Second, no edge of $\Omega_q$ can ever be extracted. Indeed, initially $\Omega_q\subseteq T(\mathscr{T})$. Assume that $\Omega_q$ is contained in a current component $S$ of the sequential algorithm. For every $l\in\Omega_q$, Lemma~\ref{lem4.1} yields
\[
m(l|_S)\le m(l|_{\Omega_q})<d+1,
\]
so $l$ is not eligible for extraction. Therefore $\Omega_q$ remains contained in every subsequent component and, in particular, $\Omega_q\subseteq C$.

The first part gives $C\subseteq\Omega_q$, whereas the second gives $\Omega_q\subseteq C$. Thus $C=\Omega_q$ and $D=D_1\cup\cdots\cup D_q$, independently of all sequential choices. The complete partition is therefore unique.

$\Box$
\end{proof}

The uniqueness of the complete partition directly yields the following four conclusions.

\begin{corollary}
    The CNDC of a T-mesh is uniquely determined.
\end{corollary}

\begin{corollary}\label{cor4.2}
    \textcolor{black}{Given a T-mesh $\mathscr{T}$, let $t$ be the number of T $l$-edges of $T(\mathscr{T})$, and let $M(T(\mathscr{T}))$ be its conformality matrix. Let $\{\overline{T}_1,\overline{T}_2\}$ be the complete partition, where the CNDC is $\overline{T}_1=\{l_1,l_2,\ldots,l_s\}$ and $s\le t$. Then}
    \begin{equation}\label{dimformula}
        \mathrm{rank}\left(M\left(T(\mathscr{T})\right)\right)=(t-s)(d+1)+\mathrm{rank}(M(\overline{T}_1)),
    \end{equation}
    where $M(\overline{T}_1)$ is the conformality matrix of the CNDC.
\end{corollary}
\begin{proof}
    Since $\overline{T}_2$ is the diagonalizable component of $T(\mathscr{T})$, its corresponding conformality matrix $M(\overline{T}_2)$ has full row rank with $\mathrm{rank}(M(\overline{T}_2)) = (t - s)(d + 1)$. Due to the block structure of the complete partition, the rank of the total matrix is the sum of the ranks of its sub-components, namely
    $$\mathrm{rank}(M(T(\mathscr{T}))) = \mathrm{rank}(M(\overline{T}_1)) + \mathrm{rank}(M(\overline{T}_2)).$$
    Substituting the rank of $M(\overline{T}_2)$ directly yields the formula.

    $\Box$
\end{proof}

The following corollary corresponds to the main theorem (Theorem 4.1) in \cite{huang2023}. However, as discussed previously, it is precisely the uniqueness of the CNDC established in Theorem~\ref{thmcompuni} that now allows this spline dimension formula to be formally and rigorously stated.

\begin{corollary}\label{cor4.3}
    \textcolor{black}{Given a T-mesh $\mathscr{T}$, let $t$ be the number of T $l$-edges of $T(\mathscr{T})$, and let $M(T(\mathscr{T}))$ be its conformality matrix. Let $\{\overline{T}_1,\overline{T}_2\}$ be the complete partition, where the CNDC is $\overline{T}_1=\{l_1,l_2,\ldots,l_s\}$ and $s\le t$. Then}
    \begin{equation}\label{dimformulanew}
        \dim \mathbf{S}_{d}(\mathscr{T})=(d+1)^2+(c+s-t)(d+1)+n_v-\mathrm{rank}(M(\overline{T}_1)),
    \end{equation}
    where $c$ is the number of all cross-cut edges of $\mathscr{T}$, $n_v$ is the number of all interior vertices of $\mathscr{T}$, and $M(\overline{T}_1)$ is the conformality matrix of the CNDC of $T(\mathscr{T})$.
\end{corollary}

By combining Theorem~\ref{thmcompuni} and Theorem~\ref{thm dimension stable}, we obtain the following necessary and sufficient condition for dimensional stability. First, the uniqueness of the CNDC guarantees that $s$, namely the number of edges within the CNDC defined in Corollary~\ref{cor4.3}, remains an invariant under structurally isomorphic mappings. Since all other topological quantities in the dimension formula \eqref{dimformulanew} are likewise invariant within the equivalence class $[\mathscr{T}]_{\mathrm I}$, the dimensional stability of the spline space depends exclusively on the CNDC.

\begin{corollary}\label{cor4.4}
    Let $[\mathscr{T}]_{\mathrm I}\in\mathrm{SIC}$. The spline space $\mathbf{S}_{d}(\mathscr{T})$ is dimensionally stable if and only if for any two T-meshes $\mathscr{T}_1, \mathscr{T}_2 \in [\mathscr{T}]_{\mathrm I}$, the ranks of the conformality matrices of their corresponding CNDCs are equal, i.e.,
    $$\textcolor{black}{\mathrm{rank}\bigl(M(\overline{T}_{\mathrm{CNDC}}^{(1)})\bigr) = \mathrm{rank}\bigl(M(\overline{T}_{\mathrm{CNDC}}^{(2)})\bigr),}$$
    where \textcolor{black}{$\overline{T}_{\mathrm{CNDC}}^{(1)}$ and $\overline{T}_{\mathrm{CNDC}}^{(2)}$} are the unique CNDCs of $\mathscr{T}_1$ and $\mathscr{T}_2$, respectively.
\end{corollary}
\begin{proof}
    According to Theorem~\ref{thm dimension stable}, the spline space is dimensionally stable if and only if $\dim \mathbf{S}_{d}(\mathscr{T})$ is constant for all meshes within the equivalence class $[\mathscr{T}]_{\mathrm I}$. In the dimension formula \eqref{dimformulanew}, the parameters $c$, $t$, and $n_v$ are identical for any meshes in $[\mathscr{T}]_{\mathrm I}$ due to their shared topology. Furthermore, Theorem~\ref{thmcompuni} guarantees that the CNDC is uniquely determined, which implies that the number of remaining edges $s$ is also a topological invariant. Therefore, the dimension remains constant if and only if \textcolor{black}{the rank of the conformality matrix of the corresponding CNDC is invariant} for all meshes in $[\mathscr{T}]_{\mathrm I}$.
    $\Box$
\end{proof}

\subsection{Characteristics of the Diagonalizable T-meshes}

In this subsection, we provide several algebraic and topological characterizations of diagonalizable T-meshes. We demonstrate that the diagonalizability of a T-mesh can be determined entirely by the rank and block structure of the conformality matrices of its partition components.

\begin{proposition}\label{prop diagonalizable}
\color{black}
A T-mesh $\mathscr{T}$ is diagonalizable if and only if it admits a $t$-partition $\{\overline l_1,\ldots,\overline l_t\}$ such that every local conformality matrix $M(\overline l_i)$ has full row rank, where $t$ is the number of T $l$-edges of $T(\mathscr{T})$. In that case,
\[
\operatorname{rank}M(T(\mathscr{T}))
=
\sum_{i=1}^{t}\operatorname{rank}M(\overline l_i)
=t(d+1).
\]
\end{proposition}
\begin{proof}
\color{black}
For a fixed $t$-partition, the matrix $M(\overline l_i)$ is a Vandermonde-type matrix with $d+1$ rows and $n(\overline l_i)$ columns. Since the retained vertex coordinates are distinct, it has full row rank exactly when $n(\overline l_i)\ge d+1$. Consequently, the full-row-rank condition for every $i$ is equivalent to the defining condition for a diagonalizable T-mesh. Under the same ordering, the global conformality matrix is block triangular with diagonal blocks $M(\overline l_i)$, so its rank is the sum of the block ranks.

$\Box$
\end{proof}

\begin{corollary}
\color{black}
A T-mesh $\mathscr{T}$ is diagonalizable if and only if there exists a $k$-partition $\{\overline T_1,\ldots,\overline T_k\}$ such that, for each component $\overline T_i$, its T $l$-edges can be ordered as
\[
l_1^{(i)}\succ\cdots\succ l_{t_i}^{(i)}
\]
and the retained edge $\overline l_j^{(i)}$, obtained after deleting the intersection vertices shared with $l_1^{(i)},\ldots,l_{j-1}^{(i)}$, satisfies
\[
n(\overline l_j^{(i)})\ge d+1,
\qquad j=1,\ldots,t_i.
\]
\end{corollary}
\begin{proof}
\color{black}
A diagonalizing order of all T $l$-edges can be divided into consecutive blocks, giving a $k$-partition with the stated property. Conversely, concatenating the admissible orders of the partition components in their block order gives a global order whose retained edges all have at least $d+1$ vertices. This is precisely a diagonalizing order.

$\Box$
\end{proof}

\begin{corollary}
\color{black}
Let $\mathscr{T}$ be diagonalizable, and let $\{\overline T_1,\ldots,\overline T_k\}$ be a $k$-partition satisfying the conditions of the preceding corollary. Then
\[
\operatorname{rank}M(T(\mathscr{T}))
=
\sum_{i=1}^{k}\operatorname{rank}M(\overline T_i).
\]
\end{corollary}

\textcolor{black}{For a $k$-partition $\{\overline{T}_1,\ldots,\overline{T}_k\}$ of $T(\mathscr{T})$, suppose that the corresponding conformality vector space decomposes as the algebraic direct sum}
\[
\mathsf{CVS}[T(\mathscr{T})]
=
\bigoplus_{i=1}^k \mathsf{CVS}[\overline{T}_i].
\]
\textcolor{black}{We then use the symbol $\biguplus$ to record this structural decomposition:}
\[
T(\mathscr{T}) = \biguplus_{i=1}^k \overline{T}_i.
\]
\textcolor{black}{For a general T-mesh $\mathscr{T}$, let $t$ be the number of T $l$-edges of $T(\mathscr{T})$, and suppose that the CNDC $\overline{T}_{\mathrm{CNDC}}$ contains $s$ of them. If $\{\overline l_1,\ldots,\overline l_{t-s}\}$ is a $(t-s)$-partition of the diagonalizable component $\overline{T}_{\mathrm D}$, then}
\[
 \textcolor{black}{T(\mathscr{T}) = \overline{T}_{\mathrm{CNDC}} \biguplus \left( \biguplus_{i=1}^{t-s} \overline{l}_i \right).}
\]

From the algebraic perspective of basis function construction, this structural decomposition explicitly governs the distribution of the spline degrees of freedom (DOFs) across the sub-structures. Taking the dimensions of the associated conformality vector spaces yields the exact additive relation:
\[
 \dim \mathsf{CVS}[T(\mathscr{T})] = \dim \mathsf{CVS}[\overline{T}_{\mathrm{CNDC}}] + \sum_{i=1}^{\textcolor{black}{t-s}} \dim \mathsf{CVS}[\overline{l}_i],
\]
where each decoupled T $l$-edge independently contributes exactly $\dim \mathsf{CVS}[\overline{l}_i] = n(\overline{l}_i)-(d+1)$ local DOFs to the global space. Consequently, the total basis set $\mathcal{B}$ spanning the spline space can be decomposed into
$$\mathcal{B} = \mathcal{B}_{\mathrm{CNDC}} \cup \left( \bigcup_{i=1}^{\textcolor{black}{t-s}} \mathcal{B}_{\overline{l}_i} \right),$$
where the cardinality of each edge-associated basis subset strictly satisfies $|\mathcal{B}_{\overline{l}_i}| = n(\overline{l}_i)-(d+1)$. This algebraic requirement mandates that each $\overline{l}_i$ requires exactly $n(\overline{l}_i)-(d+1)$ local tensor product B-splines as the basis functions corresponding to its independent DOFs. However, due to the topological confinement of surrounding T-junctions, the required local tensor product supports for these B-splines often cannot be fully accommodated within the unextended boundaries of $\mathscr{T}$. This intrinsic mismatch between the algebraic demand for $n(\overline{l}_i)-(d+1)$ basis functions and the local geometric capacity directly motivates the ``edge extension'' mechanism proposed in~\cite{zhong2026}.

\begin{example}\label{exm:decoupling_dof}
    To explicitly demonstrate this algebraic distribution of DOFs via the decomposition formula, we consider the completely diagonalizable T-mesh $\mathscr{T}$ from Fig.~\ref{fig2.sub.1} under the spline space $\textcolor{black}{\mathbf{S}}_1(\mathscr{T})$. Since $\overline{T}_{\mathrm{CNDC}} = \emptyset$, the structural decomposition yields $T(\mathscr{T}) = \biguplus_{i=1}^t \overline{l}_i$, which induces the following dimension formula:
    \[
     \dim \mathsf{CVS}[T(\mathscr{T})] = \sum_{i=1}^t \dim \mathsf{CVS}[\overline{l}_i].
    \]
    This equation dictates that every single T $l$-edge $\overline{l}_i$ must be allocated exactly $|\mathcal{B}_{\overline{l}_i}| = n(\overline{l}_i)-(d+1) = 3-2=1$ localized DOF, which requires exactly one local tensor product B-spline. Nevertheless, by verifying the local configuration of $\mathscr{T}$, certain interior T $l$-edges lack a valid local tensor product support inside the unextended boundaries of $\mathscr{T}$ to define their corresponding local B-spline basis functions. Hence, to geometrically realize the algebraic allocation $|\mathcal{B}_{\overline{l}_i}| = 1$, the mesh domain must be embedded into an expanded state via edge extensions to host the required knot spans.
\end{example}

\begin{remark}\label{rem:zhong_foundation}
    We emphasize that the direct-sum decomposition
$$  \mathsf{CVS}[T(\mathscr{T})] = \mathsf{CVS}[\overline{T}_{\mathrm{CNDC}}] \oplus \left(\bigoplus_{i=1}^{\textcolor{black}{t-s}} \mathsf{CVS}[\overline{l}_i]\right) $$
and its \textcolor{black}{dimension-additivity property} established in this section serve as the indispensable mathematical foundation for the framework in~\cite{zhong2026}. The core pipeline in~\cite{zhong2026} directly capitalizes on this result by executing sequential edge extensions $\mathscr{T} \hookrightarrow \mathrm{ext}(\mathscr{T})$ to transform any arbitrary T-mesh into a diagonalizable extended T-mesh $\mathrm{ext}(\mathscr{T}) = \biguplus_{j=1}^m \overline{l}_j$, where each extended edge is locally guaranteed a valid tensor grid hosting exactly $n(\overline{l}_j)-(d+1)$ \textcolor{black}{basis functions}. The valid spline bases over the original unextended domain $\mathscr{T}$ are subsequently recovered by computing the nullspace of the Extended Edge Elimination (EEE) conditions, formulated as a linear system. Consequently, the framework in~\cite{zhong2026} is fundamentally built upon the theoretical results established in this paper.
\end{remark}

To further clarify the relationship between our results and existing adaptive spline technologies, we present a formal comparison with the widely used LR-meshes introduced in~\cite{lr}. It is worth noting that while LR-meshes in~\cite{lr} can be defined for spline spaces with arbitrary smoothness, they are considered under the highest-order smoothness condition in this paper to align with our focus. An LR-mesh is constructed via a sequence of meshline insertions, where each step requires inserting a line segment that contains a sufficient number of vertices to support new local B-splines.

\begin{definition}[LR-mesh(\cite{lr})]\label{def:lrmesh}
    A T-mesh $\mathscr{T}$ is called an LR-mesh if it can be obtained from an initial tensor product mesh $\mathscr{T}_0$ through a finite sequence of meshline insertions $\mathscr{T}_0 \subset \mathscr{T}_1 \subset \cdots \subset \mathscr{T}_N = \mathscr{T}$. Each inserted meshline $E_k$ ($k=1,\ldots,N$) must split a set of existing elements and contain at least $d+2$ \textcolor{black}{vertices created by intersections with orthogonal meshlines or boundaries} to host at least one new local tensor product B-spline, meaning its interior vertex count satisfies $n(E_k) \ge d+2$.
\end{definition}

\begin{theorem}\label{thm:lr_is_diagonalizable}
    Every LR-mesh $\mathscr{T}$ is a diagonalizable T-mesh and it admits a full $t$-partition:
    \[
     T(\mathscr{T}) = \biguplus_{i=1}^t \overline{l}_i.
    \]
\end{theorem}

\begin{proof}
    Let $\mathscr{T}$ be an LR-mesh generated by a sequence of meshline insertions $\mathcal{E} = \{E_1, E_2, \ldots, E_N\}$. Since rays do not affect the diagonalizability of a T-mesh, the characterization of diagonalizable T-meshes depends solely on T $l$-edges. Therefore, without loss of generality, each inserted meshline $E_k$ can be identified directly as a single T $l$-edge. We then define the diagonalization sequence $\sigma$ in the exact reverse order of the insertion process, i.e.,
    \[
     \sigma = (\overline{l}_1, \overline{l}_2, \ldots, \overline{l}_N) = (E_N, E_{N-1}, \ldots, E_1).
    \]
    For each step $k = N, N-1, \ldots, 1$, since $E_k$ is the latest inserted line in the sub-structure $\mathscr{T}_k = \mathscr{T}_0 \cup \left(\bigcup_{j=1}^k E_j\right)$, no subsequent lines can terminate on $E_k$ within $\mathscr{T}_k$. Thus, $E_k$ naturally satisfies the condition of a decoupled T $l$-edge in $\mathscr{T}_k$. By induction, the reverse sequence $\sigma$ successfully eliminates all internal constraints step-by-step until reducing to $\mathscr{T}_0$, proving that $\overline{T}_{\mathrm{CNDC}} = \emptyset$ and the mesh is diagonalizable.

    $\Box$
\end{proof}

Crucially, the converse of Theorem~\ref{thm:lr_is_diagonalizable} does not hold. Diagonalizable T-meshes form a strictly larger class than LR-meshes because our direct-sum decomposition can accommodate topological structures with cyclic dependencies (T-cycles) that can never be generated via sequential line insertions. We demonstrate this fundamental distinction through the following example corresponding to the configuration in \textcolor{black}{Fig.}~\ref{fig:tcycle_mesh}.

\begin{example}\label{exm:tcycle_counterexample}
    Consider the T-mesh $\mathscr{T}$ illustrated in Fig.~\ref{fig:tcycle_mesh} under the spline space $\mathbf{S}_2(\mathscr{T})$.

    \begin{figure}[htbp]
        \centering
        \begin{tikzpicture}[scale=0.65, thick]
            \draw[black] (0,0) rectangle (8,8);

            \draw[black] (0,1) -- (8,1);
            \draw[black] (0,2) -- (8,2);
            \draw[black] (0,6) -- (8,6);
            \draw[black] (0,7) -- (8,7);

            \draw[black] (1,0) -- (1,8);
            \draw[black] (2,0) -- (2,8);
            \draw[black] (6,0) -- (6,8);
            \draw[black] (7,0) -- (7,8);

            \draw[black, ultra thick] (3,5) -- (7,5); 
            \draw[black, ultra thick] (3,3) -- (3,7); 
            \draw[black, ultra thick] (1,3) -- (5,3); 
            \draw[black, ultra thick] (5,1) -- (5,5); 

            \draw[black] (0,4) -- (3,4);

            \filldraw[black] (3,5) circle (2.5pt) node[above left] {$v_1$};
            \filldraw[black] (7,5) circle (2.5pt) node[above right] {$v_2$};
            \filldraw[black] (3,3) circle (2.5pt) node[below left] {$v_3$};
            \filldraw[black] (3,7) circle (2.5pt) node[above left] {$v_4$};
            \filldraw[black] (1,3) circle (2.5pt) node[below left] {$v_5$};
            \filldraw[black] (5,3) circle (2.5pt) node[below right] {$v_6$};
            \filldraw[black] (5,1) circle (2.5pt) node[below right] {$v_7$};
            \filldraw[black] (5,5) circle (2.5pt) node[above right] {$v_8$};

        \end{tikzpicture}
        \caption{A diagonalizable T-mesh containing a T-cycle, which is not an LR-mesh.}
        \label{fig:tcycle_mesh}
    \end{figure}

    As marked on the grid, the endpoints of the four main interior T $l$-edges are explicitly given by:
    \[
    \overline{l}_1 = \overline{v_1v_2}, \quad \overline{l}_2 = \overline{v_3v_4}, \quad \overline{l}_3 = \overline{v_5v_6}, \quad \overline{l}_4 = \overline{v_7v_8}.
    \]
    This T-mesh admits a valid diagonalization sequence $\sigma = (\overline{l}_2, \overline{l}_1, \overline{l}_4, \overline{l}_3)$. \textcolor{black}{Each retained T $l$-edge $\overline{l}_i$ satisfies $n(\overline{l}_i)\ge d+1=3$; hence the mesh is diagonalizable.}

    However, the mesh contains a clear T-cycle topological structure among the edges $\{\overline{l}_1, \overline{l}_2, \overline{l}_3, \overline{l}_4\}$ due to their mutual geometric terminations ($v_1 \in \overline{l}_2$, $v_3 \in \overline{l}_3$, $v_6 \in \overline{l}_4$, and $v_8 \in \overline{l}_1$). As a result, this configuration cannot be viewed as being formed by inserting one independent line segment at a time under Definition~\ref{def:lrmesh}, meaning $\mathscr{T}$ is \textbf{not an LR-mesh}. This confirms that the class of diagonalizable T-meshes is strictly more general than that of LR-meshes.
\end{example}

Finally, we establish the main theorem of this section, demonstrating that diagonalizability is a sufficient condition for a T-mesh to be dimensionally absolutely stable.

\begin{theorem}\label{thm:diag_absolute_stable}
\color{black}
Every diagonalizable T-mesh $\mathscr{T}$ is dimensionally absolutely stable.
\end{theorem}
\begin{proof}
\color{black}
Let $\mathrm{GT}_1,\mathrm{GT}_2\in[T(\mathscr{T})]_{\mathrm S}$, and let
\[
\beta:L(\mathrm{GT}_1)\longrightarrow L(\mathrm{GT}_2)
\]
be a structurally similar map. Choose a diagonalizing order $l_1\succ\cdots\succ l_t$ for $T(\mathscr{T})$. Structural similarity transfers this order to each $\mathrm{GT}_q$ ($q=1,2$). Because $\beta$ preserves both the number of vertices on every corresponding $l$-edge and all pairwise incidences, the same number of previously shared intersection vertices is deleted at each step. Hence the corresponding retained edges satisfy
\[
n(\overline l_i^{(1)})=n(\overline l_i^{(2)})\ge d+1,
\qquad i=1,\ldots,t.
\]
Thus both generalized components are diagonalizable in the analogous sense. Their conformality vector spaces decompose into the local edge spaces, and
\[
\dim\mathsf{CVS}[\mathrm{GT}_q]
=
\sum_{i=1}^{t}\bigl(n(\overline l_i^{(q)})-(d+1)\bigr),
\qquad q=1,2.
\]
The equality of the corresponding retained vertex counts gives
\[
\dim\mathsf{CVS}[\mathrm{GT}_1]
=
\dim\mathsf{CVS}[\mathrm{GT}_2].
\]
Therefore $\mathscr{T}$ is dimensionally absolutely stable.

$\Box$
\end{proof}

At the end of this subsection, we propose a conjecture regarding the relationship between \textcolor{black}{absolute stability} and diagonalizable T-meshes. This conjecture is significant because diagonalizable T-meshes represent the broadest class of T-meshes with respect to absolute stability, offering an initial framework for classifying dimensional stability.

\begin{conjecture}\label{conj 2}
    \textcolor{black}{A T-mesh $\mathscr{T}$ is dimensionally absolutely stable if and only if it is diagonalizable.}

\begin{center}
\begin{minipage}{0.88\linewidth}
\color{black}
Accordingly, T-meshes are divided into the dimensionally unstable class and the dimensionally stable class. The latter is further divided into
\begin{itemize}
    \item dimensionally absolutely stable T-meshes, conjecturally exactly the diagonalizable T-meshes;
    \item dimensionally stable but not absolutely stable T-meshes.
\end{itemize}
\end{minipage}
\end{center}
\end{conjecture}

\subsection{Summary of Section 4}
This section focuses on the decomposition of T-connected components. Initially, we introduce the general $k$-partition, then prove the existence and uniqueness of the complete partition proposed in \cite{huang2023}, successfully reducing the T-mesh dimensional stability problem to the rank stability of the conformality matrix associated with the CNDC. Finally, the conclusions are extended to dimensional absolute stability through structurally similar maps.

For diagonalizable T-meshes, we show that the T-connected component decomposes into independent T $l$-edges, where the global conformality vector space is the direct sum of individual local spaces. This provides a clear foundation for constructing spline basis functions. For arbitrary T-meshes, the complete partition isolates the diagonalizable part from the CNDC. Consequently, basis function construction over arbitrary T-meshes requires treating the CNDC as a single, unified entity.

In summary, this section provides a comprehensive analysis of T-mesh decompositions, laying the groundwork for analyzing dimensional stability and guiding the construction of basis functions for spline spaces.

Finally, the new symbols introduced in this section are listed in Table~\ref{tab:symbols section 4}.

\begin{table}[H]
\centering
\caption{List of Symbols in Section~\ref{sec: complete partition}\label{tab:symbols section 4}}
\begin{tabular}{|c|l|}
\hline
\textbf{Symbol} & \textbf{Explanation} \\
\hline
$\overline{T}_i$ & The $i$-th component of the $k$-partition \\
\hline
$l|_S$ & The $l$-edge $l$ together with the vertices retained in the component $S$ \\
\hline
$m(l|_S)$ & The number of relative mono-vertices of $l$ on $S$ \\
\hline
$\overline{T}_{\mathrm{CNDC}}$ & The completely non-diagonalizable component \\
\hline
$\overline{T}_{\mathrm{D}}$ & The diagonalizable part \\
\hline
$T(\mathscr{T}) = \biguplus_{i=1}^k \overline{T}_i$ & Decomposition satisfying $\mathsf{CVS}[T(\mathscr{T})] = \bigoplus_{i=1}^k \mathsf{CVS}[\overline{T}_i]$ \\
\hline
\end{tabular}
\end{table}

\section{Conclusion}\label{sec:conclusion}

This paper provides a systematic investigation into the dimensional stability of spline spaces over T-meshes, focusing on the properties and decompositions of T-connected components. Section 3 introduces the structurally isomorphic map, defining dimensional stability as an invariant within the structurally isomorphic class, contingent on the rank stability of the conformality matrix of $T(\mathscr{T})$. It further proposes absolute stability via structurally similar maps, which is independent of relative geometric positions and purely determined by the topology of the T-mesh, culminating in a classification of T-mesh stability types, as illustrated in~\eqref{eq:Tmesh-classification}. This framework enables a preliminary classification of T-meshes into unstable, stable (including \textcolor{black}{absolutely stable and stable but not absolutely stable}) categories.

Section 4 explores the decomposition of T-connected components, introducing the $k$-partition. The proof of the existence and uniqueness of the complete partition, as proposed in \cite{huang2023}, shifts the focus of T-mesh dimensional stability to the rank stability of the conformality matrix associated with the CNDC. For diagonalizable T-meshes, we establish that the T-connected component can be decomposed into a set of independent T $l$-edges, where the global conformality vector space is the direct sum of the local conformality vector spaces of each T $l$-edge in the $t$-partition. This decomposition offers a clear foundation for constructing basis functions of the spline space over diagonalizable T-meshes.

Future research on dimensional stability and related aspects should prioritize the following directions:
\begin{itemize}
    \item Prove or disprove Conjecture~\ref{conj 2} by investigating the necessary and sufficient conditions for absolute stability, thereby establishing a preliminary classification of dimensional stability for T-meshes.
    \item Clarify that the ultimate goal of dimensional stability classification is to determine an appropriate $k$-partition of the T-connected components of a given T-mesh $\mathscr{T}$ as $T(\mathscr{T}) = \biguplus_{i=1}^k \overline{T}_i$, with the primary challenge being to provide a topological characterization for each component $\overline{T}_i$ (where $1 \le i \le k$) that ensures the dimensional stability of $\mathbf{S}_{d}(\mathscr{T})$.
    \item Develop corresponding spline space basis functions with stable dimensions based on dimensional stability research, such as constructing local tensor product basis functions for diagonalizable T-meshes using the T $l$-edges within each $t$-partition, while treating the CNDC as a unified entity for arbitrary T-meshes.
\end{itemize}
These advancements will enhance the theoretical framework and practical applicability of T-mesh-based spline spaces.

\section*{Acknowledgements}
This work was supported by the Key Project of the National Natural Science Foundation of China (Nos.~12494550 and 12494555), the National Natural Science Foundation of China (No.~12371383), and the Fundamental Research Funds for the Central Universities (No.~WK0010000096). The authors declare that there are no conflicts of interest.

\bibliography{GMP2025template}

\end{document}